\documentclass[12pt]{article}
\usepackage[utf8]{inputenc}
\usepackage{amsmath}
\usepackage{amsthm}
\newtheorem{proposition}{Proposition}[section]
\newtheorem{lemma}[proposition]{Lemma}
\newtheorem{example}[proposition]{Example}
\newtheorem{definition}[proposition]{Definition}
\newtheorem{theorem}{Theorem}

\newtheorem{remark}[proposition]{Remark}

\usepackage{amssymb}
\newcommand {\R} {{\mathbb R}}

\newcommand {\N} {{\mathbb N}}

\newcommand {\ve}{{\varepsilon}}
\newcommand {\diam}{{\rm diam\,}}
\def\kasten{\hfill\null\nobreak\hfill \hbox{\vrule\vbox{\hrule width
6pt\vskip6pt\hrule}\vrule} \par\smallskip}

\title{Curvature Varifolds with Orthogonal Boundary}
\author{Ernst Kuwert  and Marius Müller}
\date{\today}

\begin{document}

\maketitle

\begin{abstract}
We consider the class ${\bf S}^m_\perp(\Omega)$ of $m$-dimensional surfaces in 
$\overline{\Omega} \subset \R^n$ which intersect $S = \partial \Omega$ orthogonally 
along the boundary. A piece of an affine $m$-plane in ${\bf S}^m_\perp(\Omega)$
is called an orthogonal slice. We prove estimates for the area by the 
$L^p$ integral of the second fundamental form in three cases: first 
when $\Omega$ admits no orthogonal slices, 
second for $m = p = 2$ if all orthogonal slices are topological disks,
and finally for all $\Omega$ if the surfaces are confined to a 
neighborhood of $S$. The orthogonality constraint has a 
weak formulation for curvature varifolds. We classify those
varifolds of vanishing curvature. As an application,
we prove for any $\Omega$ the existence of an orthogonal 
$2$-varifold which minimizes the $L^2$ curvature in the 
integer rectifiable class.
\end{abstract}

\begin{small}
\end{small}
%{\em Acknowledgements. } The authors would like to thank Victor Bangert
%(Freiburg) and Lei Liu (Wuhan) for helpful discussions.
% \end{small}

\section{Introduction}

For a bounded domain $\Omega \subset \R^n$ of class $C^2$, 
consider the class ${\bf S}^m_\perp(\Omega)$ of compact 
$C^1$ immersed submanifolds $\Sigma \hookrightarrow \overline{\Omega}$ 
of dimension $m \leq n-1$, which are $C^2$ in the interior and 
meet $S = \partial \Omega$ orthogonally along their boundary.
Denoting by $\nu^S$ the 
unit normal interior to $\Omega$, the condition means that
$\nu^S \in T\Sigma$ along $\partial \Sigma \hookrightarrow S$.
A natural variational problem in the class ${\bf S}^m_\perp(\Omega)$
is to minimize the $L^p$ curvature energy
$$
{\cal E}^p(\Sigma) = \int_\Sigma |A_\Sigma|^p\,d\mu_\Sigma
\quad \mbox{ where }1 \leq p < \infty.
$$
In \cite{AK16} Alessandroni and the first author constructed 
critical points of the Willmore energy subject to the 
orthogonality constraint, but only in the class of surfaces 
with small prescribed area. Here in Section 5 we address
the global problem in the case when the area is not fixed 
a priori. Area bounds depending on the curvature energy 
then play a key role, they are also of independent 
interest. It turns out that the area of surfaces in 
${\bf S}^m_\perp(\Omega)$ has an upper bound in terms of
the curvature energy, unless the domain $\Omega$ is special.

\begin{definition}[Orthogonal Slices] \label{defortho} 
Let $x_0+P \subset \R^n$ be an affine $m$-plane. Any 
component $\Delta \neq \emptyset$ of\,  
${\rm int} (\overline{\Omega} \cap (x_0+P))$ 
will be called an $m$-slice of $\Omega$. $\Delta$ is 
an orthogonal $m$-slice if it meets $S$ orthogonally
along its boundary, that is 
$$
\nu^S(x) \in P \quad \mbox{ for all }x \in \partial \Delta.
$$
\end{definition}

%$\overline{\Omega}_P = \overline{\Omega} \cap P$ and
%denote by $\partial \overline{\Omega}_P$ its boundary 
%as subset of $P$. We say that $P$ is orthogonal to 
%$\Omega$ if $P$ and $S$ intersect orthogonally along
%$\partial \overline{\Omega}_P$, that is $\nu^S(x) \in P$ 
%for all $x \in \partial \overline{\Omega}_P$. 
%\end{definition}

We have the following estimate. 

%example: symmetry planes
%In general if $P$ is orthogonal to $\Omega$, then 
%$\Omega_P: = {\rm int\,}\overline{\Omega}_P$ is a nonempty
%bounded open set of class $C^1$ with boundary 
%$C_P = \partial \overline{\Omega}_P$. Namely for any
%any $x \in C_P$, we have 
%$$
%\nabla (\mathrm{dist}_S|_P)(x) = P \nabla \mathrm{dist}_S(x) = \nu^S(x) \neq 0.
%$$
%Thus there is a neighborhood $U$ of $x$, such that $\Omega \cap P$ 
%is locally a domain $U^{+}$ with $C^1$ boundary near $x$, which
%implies
%$$
%\overline{\Omega}_P \cap U = \overline{U^+} \cap U \quad 
%\mbox{ and hence } \quad
%\Omega_P \cap U = U^{+}.
%$$
%The definition allows for points $x \in S \cap P$ which are 
%not in $\partial \overline{\Omega}_P$, i.e. there is a 
%neighborhood $U$ of $x$ in $P$ with $U \subset  \overline{\Omega}$. 
%In particular we must have $P \subset T_x S$, i.e. $P$ is
%tangential at $x$.

\begin{theorem}[Mass Bounds] \label{thm1} 
Assume that $\Omega \subset \R^n$ has no 
orthogonal $m$-slice. Then there is a constant 
$C < \infty$ depending on $n$, $m$, $p$ and $\Omega$, such 
that
\begin{equation}
\label{eqthm1bound}
|\Sigma| + |\partial \Sigma| \leq 
C\,\int_\Sigma |A_\Sigma|^p\,d\mu_\Sigma \quad
\mbox{ for any }\Sigma \in {\bf S}^m_\perp(\Omega).
\end{equation}
Moreover for $m = p = 2$, the Euler characteristic is also 
estimated by 
\begin{equation}
\label{eqeuler}
|\chi(\Sigma)| \leq
C\,\int_\Sigma |A_\Sigma|^2\,d\mu_\Sigma \quad
\mbox{ for any }\Sigma \in {\bf S}^2_\perp(\Omega).
\end{equation}
\end{theorem}

We also obtain bounds for general $C^2$ domains $\Omega$ if we 
restrict to surfaces $\Sigma$ which are supported in a 
neighborhood of $S=\partial \Omega$. More precisely let 
$\varrho_S > 0$ be the supremum of all $\varrho > 0$ such that 
the map
$$
S \times [0,\rho) \rightarrow U_\rho^+(S),\,
(x,r) \mapsto x +r \nu^S(x), 
$$ 
is a diffeomorphism. Here $U_\rho^+(S)$ denotes the one-sided
tubular neighborhood of $S$. 

\begin{theorem}[Bounds near the Boundary] \label{thmboundsnearboundary}
Let $\Omega \subset \R^n$ be a bounded domain of class $C^2$ with
boundary $S = \partial \Omega$. Then for any $\delta \in (0,\varrho_S)$
there is a constant $C = C(n,m,p,\Omega,\delta) < \infty$ such that
for $\Sigma \in S^m_\perp(\Omega)$ 
\begin{equation}
\label{eqboundsnearboundary}
|\Sigma| + |\partial \Sigma| \leq C \int_{\Sigma} |A_\Sigma|^p\,d\mu_\Sigma
\quad \mbox{ whenever }\Sigma \subset \overline{U^{+}_{\delta}(S)}.
\end{equation}
For $m = p = 2$ the Euler characteristic $\chi(\Sigma)$ has a corresponding bound. 
\end{theorem}

A variant of Theorem \ref{thm1} in the case $m = p = 2$ is as follows.

\begin{theorem} [Bounds for Surfaces] \label{thmdisktype} 
Let $\Omega \subset \R^n$ be a bounded domain of class $C^2$.
Assume that all orthogonal $2$-slices of $\Omega$ 
are topological disks. Then there is a constant $C < \infty$
depending on $n$ and $\Omega$ such that
\begin{equation}
|\Sigma| + |\partial \Sigma| + |\chi(\Sigma)| \leq
C\,\Big(\int_\Sigma |A_\Sigma|^2\,d\mu_\Sigma + \chi(\Sigma)^{+}\Big) \quad
\mbox{ for any }\Sigma \in {\bf S}^2_\perp(\Omega).
\end{equation}
\end{theorem}

For $\Sigma$ connected we have $\chi(\Sigma)^{+} = 1$
if $\Sigma$ is a disk and $\chi(\Sigma)^{+} = 0$ otherwise.
We note that Theorem \ref{thmdisktype} applies to any 
convex $C^2$ domain. In Chapter 4 we give an example 
showing that the assumptions in both Theorems 
\ref{thm1} and \ref{thmdisktype} cannot be dropped.\\
\\
There are a number of papers dealing with curvature-mass 
relations for confined curves and surfaces. M\"uller 
and R\"oger consider closed $2$-dimensional surfaces
with the type of a sphere in a $3$-dimensional domain, 
e.g.\,a ball \cite{MR14}. We also refer to the work 
of Pozzetta \cite{Poz20} and the references therein. 
In a different direction, Colding and Minicozzi
proved compactness results for sequences of embedded
minimal surfaces in $3$-manifolds, see for instance 
\cite{CM99}. A problem which is more closely related 
to the present paper is studied by Mondino \cite{Mon14},
proving mass bounds for $m$-varifolds in a Riemannian
manifold in terms of their curvature. He applies
the concept of curvature varifolds due to Hutchinson
\cite{Hut86}. The mass bounds are obtained under the 
condition that there exist no varifolds with zero 
curvature. In the special case $m = 2$, it is actually
sufficient to rule out smooth totally geodesic immersions,
see \cite{BK21}. In our case, we need to introduce
a weak formulation of orthogonality along the
boundary for general curvature varifolds. For this
we build on the notion of varifolds with boundary, 
defined and studied by Mantegazza \cite{Man96}.

\begin{definition} [orthogonal boundary]\label{defadmissible}
Let $V$ be an $m$-varifold in $\R^n$ with curvature $B \in L^1(V)$
and support in $\overline{\Omega}$, and let $\Gamma$ be an
$(m-1)$-varifold in $S = \partial \Omega$. We say that $V$ is
orthogonal to $S$ along $\Gamma$ if for all
$\phi \in C^1(\R^n \times \R^{n \times n},\R^n)$
\begin{eqnarray}
\label{eqadmissible}
&& \int_{\R^n \times G(m,n)} \big(D_P\phi \cdot B
+ \langle {\rm tr\,}B,\phi \rangle + \langle D_x\phi,P \rangle \big)\,dV(x,P)\\
\nonumber
& = &
- \int_{G_{m-1}(TS)} \langle \nu^S(x), \phi(x,\nu^S(x) \wedge Q) \rangle\,d\Gamma(x,Q).
\end{eqnarray}
\end{definition}

Here $\nu^S(x) \wedge Q$ denotes the $m$-plane spanned by $\nu^S(x)$
and $Q$. The varifold $\Gamma$ in the definition is uniquely determined,
and the class of varifolds $V$ satisfying the definition (for 
some $\Gamma$) is denoted by ${\bf CV}^m_\perp(\Omega)$.\\
\\
In Section 2 we classify varifolds $V \in {\bf CV}^m_\perp(\Omega)$ 
having curvature $B_V = 0$; using disintegration we derive a 
formula to which only the slices orthogonal to $\Omega$ contribute.
In Section 3 we then consider the infimum 
\begin{equation}
\label{eqinfimum}
\varkappa^{m,p}(\Omega) = \inf\big\{\|B_V\|_{L^p(V)}^p: 
V \in {\bf CV}^m_\perp(\Omega),\,{\bf M}(V) = 1\big\}.
\end{equation}
We prove the following generalization of Theorem \ref{thm1}.

\begin{theorem} [Varifold Bounds]  \label{thmvarifolds} Let $\Omega \subset \R^n$ be 
a bounded domain of class $C^2$. If $\Omega$ has no 
orthogonal $m$-slices, then $\varkappa^{m,p}(\Omega) > 0$ 
for all $p \in [1,\infty)$. Moreover for any 
$V \in {\bf CV}^m_\perp(\Omega)$ with associated boundary 
varifold $\Gamma$, we have 
\begin{eqnarray}
\label{eqmassbound}
{\bf M}(V) & \leq & \frac{1}{\varkappa^{m,p}(\Omega)}\,\|B_V\|_{L^p(V)}^p,\\
{\bf M}(\Gamma) & \leq & C\,\|B_V\|_{L^p(V)}^p \quad 
\mbox{ where }C = C(n,m,p,\Omega).
\end{eqnarray}
\end{theorem}

In the final section this result is applied to a 
variational problem.

\begin{theorem} [Existence of Minimizers] \label{thmexistence}
For a bounded domain $\Omega \subset \R^n$ of class $C^3$, let 
${\bf ICV}^2_\perp(\Omega)$ be the class of varifolds in 
${\bf CV}^2_\perp(\Omega)$ which are integer rectifiable. 
Then the following infimum is attained: 
\begin{equation}
\label{eqminimizer}
\varkappa(\Omega) = \inf\big\{\|B_V\|_{L^2(V)}^2:
V \in {\bf ICV}^2_\perp(\Omega),\,V \neq 0\big\}.
\end{equation}
\end{theorem}

The fact that the mass is not prescribed a priori poses a
key difficulty. But since any orthogonal $2$-slice is a 
minimizer, we may restrict to the case when such slices 
do not exist. Then Theorem \ref{thmvarifolds} applies to 
give an upper mass bound. The lower bound is proved by 
contradiction. We show that if a minimizing sequence 
$V_k$ would converge to zero then
$$
\liminf_{k \to \infty} \int_{\overline{\Omega}} |B_{V_k}|^2\,dV_k
\geq 8\pi = \int_{{\mathbb S}^2_{+}} 
|B_{{\mathbb S}^2}|^2\,d\mu_{{\mathbb S}^2},
$$
where ${\mathbb S}^2_+$ denotes the upper half sphere.
However, a comparison surface strictly below $8\pi$ is
constructed in Lemma \ref{lemmacomparison}. We note that
our proof of Theorem \ref{thmexistence} would apply to 
$m$-varifolds with curvature in $L^m$, once the round
$m$-sphere was known to minimize $L^m$ curvature 
among varifolds in $\R^n$. 
%We expect that for $m = n-1$ 
%the argument involving the Gau{\ss} image generalizes 
%to curvature varifolds.

%In Section 6 we show that the property of having an 
%orthogonal $m$-slice is preserved under $C^2$ convergence
%of surfaces $S_i \to S$. 
In Section 6 we prove that any 
compact $C^2$ hypersurface $S \subset \R^n$ has a $C^2$ 
approximation by smooth surfaces $S_i$ which have no orthogonal
$m$-slices. In particular the assumption of Theorem
\ref{thm1} and Theorem \ref{thmvarifolds} is generically
satisfied.\\
\\
{\em Acknowledgements. } The authors would like to thank 
Victor Bangert (Freiburg) and Lei Liu (Wuhan) for valuable discussions. We also thank the referees for their careful reading and helpful comments.

\section{Preliminaries} 
We derive in this section
a weak formulation of the orthogonality constraint which applies 
to varifolds. For this we recall briefly the definition of 
curvature varifolds due to Hutchinson \cite{Hut86}, see also
Mantegazza \cite{Man96}.\\
\\
For an open set $U \subset \R^n$, we consider first
the class of $m$-dimensional, properly embedded submanifolds 
$\Sigma \subset U$. The Grassmannian $G(m,n)$ of $m$-dimensional
subspaces of $\R^n$ will be identified with the set of orthogonal 
projections $P \in \R^{n \times n}$ of rank $m$. On $\Sigma$ 
we then have the associated Gau{\ss} map
\begin{equation}
\label{eqgaussmap}
P:\Sigma \to G(m,n) \subset \R^{n \times n},\,P(x) = T_x\Sigma.
\end{equation}
The tangential and normal components of a vector $X \in \R^n$ 
are denoted by $PX = X^\top$ and $P^\perp X = X^\perp$. In 
the space $BL(\R^n \times \R^n,\R^n)$ of $\R^n$-valued 
bilinear forms, we define
\begin{equation}
\label{eqbform}
B(x)(v,w) = (D_{v^\top}P)(x)w \quad 
\mbox{ for }x \in \Sigma, \mbox{ and }\,v,w \in \R^n.
\end{equation}
The second fundamental form $A(x)$ is a component of $B(x)$. 
More precisely 
\begin{equation}
\label{eqsecondfundamentalform}
B(x)(v,w)^\perp = 
D_v(Pw)(x)^\perp = A(x)(v,w) \quad \mbox{ for }v,w \in T_x\Sigma.
\end{equation}
Reversely, $B(x)$ can be expressed by $A(x)$:
for a tangent frame $\tau_\alpha = \tau_\alpha(x)$, 
$1 \leq \alpha \leq m$, and fixed $v \in T_x\Sigma$, 
$w \in \R^n$ we compute 
\begin{eqnarray*}
D_v P \cdot w & = & D_v (Pw)\\
& = & [D_v (Pw)]^\perp
+ \sum_{\alpha = 1}^m \langle D_v (Pw),\tau_\alpha \rangle\, \tau_\alpha\\
& = & A(v,w^\top) 
- \sum_{\alpha = 1}^m \langle D_v (P^\perp w),\tau_\alpha \rangle\, \tau_\alpha
\quad (\mbox{ as }D_v w = 0)\\
& = & A(v,w^\top) + \sum_{\alpha =1}^m 
\langle w^\perp, D_v \tau_\alpha \rangle\, \tau_\alpha.
\end{eqnarray*}
By definition of the second fundamental form, we conclude 
%\begin{equation}
%\label{eqgaussmapderivative}
%D_\tau P \cdot w = A(\tau,w^\top) 
%+ \sum_{\alpha =1}^m \langle A(\tau,\tau_\alpha),w^\perp \rangle\, \tau_\alpha.
%\end{equation}
%In summary, for any $v,w \in \R^p$ we conclude 
\begin{equation}
\label{eqsecondff}
B(v,w) = A(v^\top,w^\top) + \sum_{\alpha=1}^m
\langle A(v^\top,\tau_\alpha),w^\perp \rangle\, \tau_\alpha
\quad \mbox{ for all }v,w \in \R^n.
\end{equation}
In particular 
$|B(v,w)|^2 = 
|A(v^\top,w^\top)|^2 + \sum_{\alpha = 1}^m \langle A(v^\top,\tau_\alpha),w \rangle^2$.
Choosing an orthonormal basis $v_k$ of $\R^n$ with $v_\alpha = \tau_\alpha$
for $\alpha = 1,\ldots,m$, we obtain 
\begin{equation}
\label{eqcurvaturesquare}
|B|^2 = \sum_{\alpha,\beta =1}^m |A(\tau_\alpha,\tau_\beta)|^2
+ \sum_{k= m+1}^n \sum_{\alpha,\beta=1}^m \langle A(\tau_\alpha,\tau_\beta),v_k \rangle^2
= 2\,|A|^2.
\end{equation}
Furthermore, the mean curvature vector is the trace of $B$, i.e. 
\begin{equation}
\label{eqtraceB}
{\rm tr\,}B = \sum_{k=1}^n B(v_k,v_k) 
= \sum_{\alpha =1}^m A(\tau_\alpha,\tau_\alpha) = H.
\end{equation}
For completeness we note also the other traces, 
$$
\sum_{k=1}^n \langle B(v_k,w),v_k \rangle =  \langle H,w \rangle
\quad \mbox{ and } \quad \sum_{k=1}^n \langle B(v,v_k),v_k\rangle = 0.
$$
To derive the first variation formula on $\Sigma$, we now compute
\begin{equation}
\label{eqtangdivergenz}
\mathrm{div}_\Sigma\, X^\top =  \mathrm{div}_\Sigma\, X + \langle H,X \rangle 
\quad \mbox{ for any }X \in C^1(\Sigma,\R^n).
\end{equation}
%\begin{eqnarray*}
%\mathrm{div}_M\, X^\top & = & \sum_{\alpha =1}^m 
%\langle D_{\tau_\alpha} (PX),\tau_\alpha \rangle\\
%& = & \sum_{\alpha =1}^m \langle D_{\tau_\alpha} X,\tau_\alpha \rangle 
%+ \sum_{\alpha=1}^m \langle X, A(\tau_\alpha,\tau_\alpha) \rangle\\
%& = & \mathrm{div}_M\, X + \langle H,X \rangle.
%\end{eqnarray*}
%F\"ur $\varphi:M \to \R$ folgt weiter
%\begin{equation}
%\label{eqtangdivergenz2}
%\mathrm{div}_M\, (\varphi X^\top) = \langle \nabla_M \varphi, X \rangle 
%+ \varphi\, \big(\mathrm{div}_M\, X + \langle H,X \rangle\big).
%\end{equation}
Substituting $X(x) = \phi(x,P(x))$ where
$\phi \in C^1(\overline{\Omega} \times \R^{n \times n},\R^n)$,
we get 
%Perhaps, defining $\phi$ only on the Grassmanian $G(m,n)$ 
%would be more appropriate, but we assume it is extended 
%to $\R^{n \times n}$. 
$$
\mathrm{div}_\Sigma\, X = 
\big \langle D_{\tau_\alpha} X, \tau_\alpha \big\rangle
= \big \langle D_x \phi(x,P) \tau_\alpha, \tau_\alpha \big\rangle
+ \big \langle D_P \phi(x,P) D_{\tau_\alpha}P,\tau_\alpha \big\rangle.
$$
Using the Hilbert-Schmidt scalar product, we can write 
$$
\big \langle D_x \phi(x,P) \tau_\alpha, \tau_\alpha \big\rangle 
%= {\rm tr}_P D_x\phi(x,P) 
= \langle D_x\phi(x,P),P \rangle.
$$
Furthermore, recalling the definiton of $B$ we have
\begin{eqnarray*}
\big \langle D_P \phi(x,P) D_{\tau_\alpha}P,\tau_\alpha \big\rangle
& = & \big \langle D_P \phi(x,P) B(\tau_\alpha,\cdot),\tau_\alpha \big \rangle\\ 
& = & \big \langle D_P \phi(x,P) B(e_i,\cdot),e_i \big \rangle\\
& = & \partial_{P^k_j}\phi^i(x,P) B_{ij}^k\\
& = : & D_P\phi(x,P) \cdot B.
\end{eqnarray*}
Collecting terms we obtain from \eqref{eqtangdivergenz}, for $X(x) = \phi(x,P(x))$,
\begin{equation}\label{eq:divsigmaXtangential}
\mathrm{div}_\Sigma X^T  = 
D_P\phi(x,P) \cdot B + \langle {\rm tr\,}B,\phi(x,P) \rangle
+ \langle D_x\phi(x,P),P \rangle. 
\end{equation}
Now associate to $\Sigma$ the $m$-varifold $V = V_\Sigma$ in $U$ 
with weight measure $\mu_\Sigma = {\cal H}^m \llcorner \Sigma$ and 
measure $V^x = \delta_{T_x\Sigma}$ on $G(m,n)$, for each $x \in \Sigma$.
If $X(x) = \phi(x,P(x))$ has compact support in $\Sigma$, then 
\begin{eqnarray*}
\int_\Sigma \mathrm{div}_\Sigma X^\top \,d\mu_\Sigma =
 \int_{U \times G(m,n)}\!\!\!\!
\big(D_P\phi(x,P) \cdot B + \langle {\rm tr\,}B,\phi(x,P) \rangle
 + \langle D_x\phi(x,P),P \rangle\big)\,dV(x,P). 
\end{eqnarray*}
On the basis of this calculation Hutchinson gave the following definition. 

\begin{definition} [varifold curvature] \label{defvarifoldcurvature}
Let $V$ be a varifold in $U \subset \R^n$. We say that $V$ has 
weak curvature $B \in L^1_{{\rm loc}}(V)$, where 
$B(x,P) \in BL(\R^n \times \R^n,\R^n)$,
if for any $\phi \in C^1(U \times \R^{n \times n},\R^n)$ 
with compact support in $U$ we have the identity
\begin{equation}
\label{eqweakcurvature}
\int_{U \times G(m,n)} \hspace{-11mm}
\big(D_P\phi(x,P) \cdot B + \langle {\rm tr\,}B,\phi(x,P) \rangle
 + \langle D_x\phi(x,P),P \rangle\big)\,dV(x,P) = 0.
\end{equation}
\end{definition}
%One can check that $B$ is unique if it exists \cite{Hut86}.
For test functions $\phi = \phi(x)$ the identity simplifies 
as follows, see Lemma 38.4 in \cite{Sim84}, 
\begin{eqnarray*}
\delta V(\phi) & = & \int_{U \times G(m,n)} \langle D\phi(x),P \rangle\,dV(x,P)\\
& = &  - \int_{U \times G(m,n)} \langle {\rm tr\,}B(x,P), \phi(x) \rangle \,dV(x,P)\\
& = & - \int_U \Big\langle \int_{G(m,n)} {\rm tr\,}B(x,P)\,dV^x(P),\, 
\phi(x) \Big \rangle\,d\mu_V(x).
\end{eqnarray*} 
As $|{\rm tr\,}B| \leq \sqrt{m} \,|B|$, we have for $K = {\rm spt\,}\phi \times G(m,n)$
$$
|\delta V(\phi)| \leq \sqrt{m} \,\|B\|_{L^1(V,K)} \,\|\phi\|_{C^0(U)}. 
$$
%\begin{eqnarray*}
%\int_U \Big|\int_{G(m,n)} {\rm tr\,}B(x,P)\,dV^x(P)\Big|\,d\mu_V(x)
%& \leq & \sqrt{m} \int_U \int_{G(m,n)} |B(x,P)|\,dV^x(P)\,d\mu_V(x)\\
%& = & \sqrt{m}\, \|B\|_{L^1(V)}.
%\end{eqnarray*}
Therefore Definition \ref{defvarifoldcurvature} implies that $V$ has 
weak mean curvature in $U$ given by 
%see \cite{Allar} 
\begin{equation}
\label{eqmeancurvatureinduced}
H_V(x) =  \int_{G(m,n)} {\rm tr\,}B(x,P)\,dV^x(P) \in 
L^1_{{\rm loc}}(\mu_V).
\end{equation}
Now assume that $\Sigma$ is a compact $C^1$-submanifold with boundary
$\partial \Sigma$, which is $C^2$ in the interior and with second fundamental form $A_\Sigma \in L^p(\mu_\Sigma)$, $p \geq 1$. Denote by $\eta$ the interior co-normal and by
$\sigma_{\partial \Sigma}$ the induced boundary measure. We approximate $\Sigma$ with compact $C^1$-subdomains $\Sigma_\varepsilon \subset {\rm int}(\Sigma)$ such that $\partial \Sigma_\epsilon$ goes to $\partial \Sigma$ in the $C^1$-topology.  
 Notice that 
$$
 \int_{\Sigma_\varepsilon} \mathrm{div}_\Sigma X^\top \,d\mu_\Sigma
= -  \int_{\partial \Sigma_\varepsilon} \langle X,\eta_\epsilon \rangle\,d\sigma_{\partial \Sigma_\varepsilon} \overset{\varepsilon \rightarrow 0}{ \longrightarrow } -  \int_{\partial \Sigma} \langle X,\eta\rangle\,d\sigma_{\partial \Sigma}.
$$
Inserting \eqref{eq:divsigmaXtangential} and passing to the limit 
$\varepsilon \rightarrow 0$, thereby using $A_\Sigma \in L^p$, we obtain  
for $V= V_\Sigma$ the formula of Mantegazza 
(see Section 3 in \cite{Man96})
\begin{equation*}
\int_{U \times G(m,n)}\!\!\!\!
\big(D_P\phi(x,P) \cdot B + \langle {\rm tr\,}B,\phi(x,P) \rangle
 + \langle D_x\phi(x,P),P \rangle\big)\,dV(x,P) = - \int_{\partial \Sigma} \langle X , \eta \rangle \, d
\sigma_{\partial \Sigma}. \end{equation*}

Our interest is in the case when $\Sigma \subset \overline{\Omega}$
meets $S = \partial \Omega$ orthogonally along the boundary. This is 
true if and only if the last equation holds true with 
$\eta = \nu^S\vert_{\partial \Sigma}$. Motivated by this we derive 
the following weak formulation. 
%$$
%\int_\Sigma \big(\mathrm{div}_\Sigma X + \langle \vec{H}, X \rangle\big)\,d\mu_\Sigma
%\int_\Sigma \mathrm{div}_\Sigma X^\top \,d\mu_\Sigma
%= - \int_{\partial \Sigma} \langle X ,\nu^S \rangle\,d\sigma_{\partial \Sigma}.
%$$
We associate to $\partial \Sigma$ the $(m-1)$-varifold $\Gamma$,
i.e. $\Gamma$ is a Radon measure on the bundle 
$$
G_{m-1}(TS) = \{(x,Q): x \in S,\,Q \in G(m-1,n),\,Q \subset T_x S\}.
$$
For $p:G_{m-1}(TS) \to S$, $p(x,Q) = x$, we have the weight measure
$$
p_\ast \Gamma = {\cal H}^{m-1} \llcorner \partial \Sigma =:\sigma_\Gamma,
$$
and the vertical measure $\Gamma^x = \delta_{T_x(\partial \Sigma)}$ 
for $x \in \partial \Sigma$. Thus for $\psi \in C^0(G_{m-1}(TS))$
%\begin{eqnarray*}
%\int_\Sigma (\mathrm{div}_\Sigma X + \langle H, X \rangle)\,\d\mu_\Sigma
%& = & - \int_\Gamma \langle X,\eta \rangle\,d{\cal H}^{m-1}_{\partial \Sigma}\\
%& = & - \int_\Gamma \langle X,\nu^S \rangle\,d{\cal H}^{m-1}_{\partial \Sigma}.
%\end{eqnarray*}
$$
\int_{\partial \Sigma} \psi(x,T_x(\partial \Sigma))\,d\sigma_{\partial \Sigma}(x) 
= \int_S \int_{G_{m-1}(T_x S)} \psi(x,Q)\,d\Gamma^x(Q)\,d\sigma_\Gamma(x) 
= \Gamma(\psi).
$$
In order to state the orthogonality condition, we adopt 
the following notation: for $Q \in G_{m-1}(\R^n)$ and $v \in \R^n$ 
with $v \notin Q$, we denote by $v \wedge Q \in G_m(\R^n)$ the 
span of $Q$ and $v$. In particular if $\Sigma$ meets $S$ orthogonally
then $P(x) = \nu^S(x) \wedge T_x(\partial \Sigma)$, and for any 
$\phi \in C^1(\overline{\Omega} \times \R^{n \times n},\R^n)$ 
we can write putting $X(x) = \phi(x,P(x))$
\begin{eqnarray*}
\int_{\partial \Sigma} \langle \nu^S,X \rangle\,d\sigma_{\partial \Sigma}
& = & 
\int_{\partial \Sigma} \big\langle \nu^S(x),
\phi\big(x,\nu^S(x) \wedge T_x(\partial \Sigma)\big) \big\rangle\,d\sigma_\Gamma(x)\\ 
& = & \int_{G_{m-1}(TS)} \langle \nu^S(x),\phi(x,\nu^S(x) \wedge Q)\rangle\,d\Gamma(x,Q).
\end{eqnarray*}

This suggests % the following notion. 

\begin{definition} [orthogonal boundary]\label{defadmissible}
Let $V$ be an $m$-varifold in $\R^n$ with curvature $B \in L^1(V)$
and support in $\overline{\Omega}$, and let $\Gamma$ be an 
$(m-1)$-varifold in $S = \partial \Omega$. We say that $V$ is 
orthogonal to $S$ along $\Gamma$ if for all 
$\phi \in C^1(\R^n \times \R^{n \times n},\R^n)$
\begin{eqnarray} 
\label{eqadmissible}
&& \int_{\R^n \times G(m,n)} \big(D_P\phi \cdot B
+ \langle {\rm tr\,}B,\phi \rangle + \langle D_x\phi,P \rangle \big)\,dV(x,P)\\
\nonumber
& = &  
- \int_{G_{m-1}(TS)} \langle \nu^S(x), \phi(x,\nu^S(x) \wedge Q) \rangle\,d\Gamma(x,Q).
\end{eqnarray} 
\end{definition}
The class of varifolds $V$ satisfying the definition (for 
some $\Gamma$) is denoted by ${\bf CV}^m_\perp(\Omega)$.
%As ${\rm spt\,}V \subset \overline{\Omega}$ is assumed,
%the first integral is actually only on $\overline{\Omega}$.
\begin{remark} \label{rmkfirstvariation} 
{\em A notion of orthogonality at the boundary 
can be defined also for varifolds $V$ which 
only have mean curvature in $L^1(\mu_V)$. For a given 
Radon measure $\sigma$ on $S$, one says that $V$ is orthogonal 
along $\sigma$ if for all $\phi \in C^1_c(\R^n,\R^n)$ 
$$
\delta V(\phi) = \int_{\R^n \times G(m,n)} \langle D\phi,P \rangle\,dV(x,P)
= - \int_{\R^n} \langle H_V,\phi \rangle\,d\mu_V
  - \int_S \langle \nu^S(x),\phi(x) \rangle\,d\sigma(x).
$$
Definition \ref{defadmissible} implies this property 
for $H_V$ is as in (\ref{eqmeancurvatureinduced}) 
and $\sigma = \sigma_{\Gamma}$, the weight measure of 
$\Gamma$. To employ test functions depending on $P$,
Definition \ref{defadmissible} will be needed.}
\end{remark}

\begin{lemma} The boundary varifold $\Gamma$ in Definition 
\ref{defadmissible} is unique if it exists.
\end{lemma}

{\em Proof. }Subtracting the identities, we get by 
Radon-Nikodym a signed measure 
$\Lambda$ on $G_{m-1}(TS)$ such that for all 
$\phi \in C^1(\R^n \times \R^{n \times n},\R^n)$
$$
\int_{G_{m-1}(TS)} \langle \nu^S(x),\phi(x,\nu^S(x) \wedge Q)\rangle\,
d\Lambda(x,Q) = 0.
$$
For given $\varphi \in C^0(G_{m-1}(TS))$ we define 
$\phi(x,P) =  \varphi(x,P \cap T_x S)\,\nu^S(x)$ whenever 
$\nu^S(x) \in P$. Using a suitable extension we obtain
a function $\phi \in C^0(\R^n \times \R^{n \times n},\R^n)$, 
such that for all $Q \in G_{m-1}(T_x S)$ 
$$
\langle \nu^S(x),\phi(x,\nu^S(x) \wedge Q)\rangle = 
\varphi(x,(\nu^S(x) \wedge Q) \cap T_x S) =
\varphi(x,Q).
$$
We conclude 
$$
\int_{G_{m-1}(TS)} \varphi(x,Q)\,d\Lambda(x,Q) = 0
\quad \mbox{ for all }\varphi \in C^0(G_{m-1}(TS)).
$$
By density this implies $\Lambda = 0$. \kasten

\begin{remark} \label{remarkscaling} 
{\em If Definition \ref{defadmissible} is valid 
with $V$ and $\Gamma$, then it holds true for $\lambda V$ and 
$\lambda \Gamma$ with the same $B$, for any $\lambda > 0$. 
However, one has the scaling property 
\begin{equation}
\label{eqscaling}
\|B\|_{L^p(\lambda V)} = \lambda^{\frac{1}{p}}\, \|B\|_{L^p(V)}.
\end{equation}
}
\end{remark}

\begin{lemma} [Comparable Masses] \label{lemmamasses} 
There exists a constant $C = C(n,\Omega) < \infty$ such 
that for any $V,\Gamma$ as in Definition \ref{defadmissible} 
one has 
\begin{eqnarray}
\label{equpperareabound}
{\bf M}(V) & \leq & C\,\big({\bf M}(\Gamma) + \|B\|_{L^1(V)}\big),\\
\label{eqlowerareabound}
{\bf M}(\Gamma) & \leq & C\,\big({\bf M}(V) + \|B\|_{L^1(V)}\big).
\end{eqnarray}
\end{lemma}

{\em Proof. }Taking $\phi(x) = x-x_0$ for $x_0 \in \Omega$ in 
Definition \ref{defadmissible} we obtain
$$
m {\bf M}(V) \leq \sqrt{m} \,\diam(\Omega)\, ({\bf M}(\Gamma) + \|B\|_{L^1(V)}).
$$
For the reverse inquality, we take $\phi = \varphi \nabla d_S$,
where $d_S$ is the signed distance from $S$ (positive in $\Omega$) and
$0 \leq \varphi \leq 1$ is a cutoff function with $\varphi|_S = 1$
and support in $U_\delta(S)$, for appropriate 
$\delta = \delta(\Omega) \leq 1$. This yields
\begin{eqnarray*}
{\bf M}(\Gamma) & = & - \int_{G_{m-1}(TS)} 
\langle \nu^S, \varphi \nabla d_S \rangle\,d\Gamma\\
& = & \int_{\overline{\Omega} \times G(m,n)} 
\big(\langle {\rm tr\,}B,\varphi \nabla d_S \rangle 
+ \langle D_x(\varphi \nabla d_S),P \rangle\big) \,dV\\
& \leq & C \big(\|B\|_{L^1(V)} + {\bf M}(V)\big).
\end{eqnarray*}
\kasten

\begin{remark}\label{rmkyoung}  {\em By Young's inequality we 
have for any $p \in (1,\infty)$
$$
\|B\|_{L^1(V)} \leq \frac{1}{p}\, \ve^{-p}\, \|B\|_{L^p(V)}^p
+ \frac{p-1}{p}\, \ve^{\frac{p}{p-1}}\, {\bf M}(V). 
$$
Therefore the inequalities (\ref{equpperareabound}),
(\ref{eqlowerareabound}) hold with $\|B\|_{L^1(V)}$ 
replaced by $\|B\|_{L^p(V)}^p$, with a constant 
that depends additionally on $p$. 
}
\end{remark}

\begin{lemma}[compactness] \label{lemmacompactness} 
Let $V_k$, $\Gamma_k$ be varifolds as in Definition 
\ref{defadmissible}, such that ${\bf M}(V_k) \leq C$ and
for some $p \in (1,\infty]$
\begin{equation}
\label{eqcurvaturebound}
\|B_k\|_{L^p(V_k)} \leq \Lambda < \infty 
\quad \mbox{ for all } k.
\end{equation}
After passing to a subsequence, one has $V_k \to V$ and 
$\Gamma_k \to \Gamma$ as varifolds in $\R^n$, and $V$, $\Gamma$ 
satisfy Definition \ref{defadmissible} for some $B \in L^p(V)$
with $\|B\|_{L^p(V)} \leq \Lambda$.
\end{lemma}

{\em Proof. }From Lemma \ref{lemmamasses} we have 
$$
{\bf M}(\Gamma_k) \leq C\,\big({\bf M}(V_k) + \|B_k\|_{L^1(V_k)}\big) 
\leq C\,\big({\bf M}(V_k) + {\bf M}(V_k)^{1-\frac{1}{p}}\, \|B_k\|_{L^p(V_k)}\big)
\leq C.
$$
Thus after passing to a subsequence, $V_k \to V$ and
$\Gamma_k \to \Gamma$ as varifolds in $\R^n$. Let 
$BL = BL(\R^n \times \R^n,\R^n)$, and introduce the functionals
$$
L_k:C^0_c(\R^n \times G(m,n),BL) \to \R,\,
L_k(\psi) = \int_{\R^n \times G(m,n)} 
\langle B_k, \psi \rangle\,dV_k.
$$
By H\"older's inequality
\begin{eqnarray*}
L_k(\psi) & \leq & \|B_k\|_{L^p(V_k)} \|\psi\|_{L^{\frac{p}{p-1}}(V_k)}\\ 
& \leq & C \|B_k\|_{L^p(V_k)} {\bf M}(V_k)^{1-\frac{1}{p}} 
\|\psi\|_{C^0(\R^n \times G(m,n))}\\
& \leq & C\,\|\psi\|_{C^0(\R^n \times G(m,n))}.
\end{eqnarray*}
After passing to a subsequence, we have $L_k \to L$ in 
$C^0_c(\overline{\Omega} \times G(m,n))'$, and 
$$
L(\psi) \leq \Lambda\,\|\psi\|_{L^{\frac{p}{p-1}}(V)}. 
$$
As $L^{\frac{p}{p-1}}(V)' = L^p(V)$ for any $p \in (1,\infty]$,
the functional $L$ is represented by some $B \in L^p(V)$.
The lemma is proved. \kasten 

An immediate consequence of Lemma \ref{lemmacompactness} is the 
following 

\begin{theorem}[Existence of Minimizers] \label{thmminimizer}
Let $\Omega \subset \R^n$ be a bounded $C^2$ domain. For
any $m \leq n-1$ and $p \in (1,\infty]$, let 
\begin{equation}
\label{eqminimizer} 
\varkappa^{m,p}(\Omega) = \inf \{\|B_V\|^p_{L^p(V)}: 
V \in {\bf CV}^m_{\perp}(\Omega), {\bf M}(V) = 1\}.
\end{equation}
Then $\varkappa^{m,p}(\Omega) < \infty$, and the infimum is attained.
\end{theorem} 

\begin{remark}{\em In view of Remark \ref{remarkscaling} we have 
\begin{equation}
\label{eqlinearbound}
\|B_V\|^p_{L^p(V)} \geq \varkappa^{m,p}(\Omega)\, {\bf M}(V)
\quad \mbox{ for all }V \in {\bf CV}^m_\perp(\Omega). 
\end{equation} 
}
\end{remark}

\begin{remark}{\em A corresponding existence result holds 
in the class ${\bf ICV}^m_\perp(\Omega)$ of integer
rectifiable varifolds in ${\bf CV}^m_\perp(\Omega)$, again for prescribed mass ${\bf M}(V)$. 
This follows by Allard's integral compactness 
theorem \cite{All72}, see also Section 5.
}
\end{remark}

\section{The Case of Zero Curvature} 

In this section we assume that $V$ is an $m$-varifold in $\R^n$
with support in $\overline{\Omega} \times G(m,n)$. We can then 
use test functions which are continuous on $\R^n$, since 
$\overline{\Omega}$ is compact. We now state from the 
appendix the following disintegration procedure for varifolds. 

\begin{itemize}
\item Let $\nu = (p_2)_\ast V$, where $p_2$ is the projection
onto $G(m,n)$. For any $P \in G(m,n)$, up to a set $E \subset G(m,n)$ 
with $\nu(E) = 0$, there exists a Radon measure $\mu_P$ on $\R^n$ with 
$\mu_P(\R^n) = 1$ such that for all $\phi \in C^0(\R^n \times G(m,n))$
one has  
$$
\int_{\R^n \times G(m,n)} \phi(x,P)\,dV(x,P) =
\int_{G(m,n)} \int_{\R^n} \phi(x,P)\,d\mu_P(x)\,d\nu(P). 
$$
\item Put $\mu_P^\perp = (P^\perp)_\ast \mu_P$ for $P \in G(m,n) \backslash E$.
Then for any $z \in P^\perp$, up to a set $E_P \subset P^\perp$ with 
$\mu_P^\perp(E_P) = 0$, there is a Radon measure $\mu_{P,z}$ on 
$z+P$ with $\mu_{P,z}(z+P) = 1$ such that for all $\varphi \in C^0(\R^n)$
$$
\int_{\R^n} \varphi(x)\,d\mu_P(x) =
\int_{P^\perp} \int_{z+P} \varphi(x)\,d\mu_{P,z}(x)\,d\mu_P^\perp(z).
$$
\end{itemize}
Taking $\phi(x,P) = \chi_E(P)$ we see  that
$V(\R^n \times E) = \nu(E) = 0$.
Likewise we compute for the set $A =\{(x,P): P \notin E,\,P^\perp x \in E_P\}$ 
\begin{eqnarray*}
V(A) & = & \int_{G(m,n)} \int_{P^\perp} 
\int_{z+P} \chi_A(x,P) \,d\mu_{P,z}(x) \,d\mu_P^\perp(z)\,d\nu(P)\\
& = &  \int_{G(m,n)} 
\underbrace{\int_{P^\perp} \chi_{E_P}(z)\,\mu_{P,z}(z+P)\,d\mu_P^\perp(z)}_{=0}
\,d\nu(P) = 0.
\end{eqnarray*}

\begin{remark}\label{rmksupports} {\em We will assume without loss of generality that the
$\mu_P$ have support in the compact set 
$K_P = \{x \in \R^n: (x,P) \in {\rm spt\,}V\} \subset \overline{\Omega}$. 
In fact 
\begin{eqnarray*}
\int_{G(m,n)} \int_{\R^n} \phi(x,P) \chi_{K_P}(x) d\mu_P(x)\,d\nu(P)
& = & \int_{\R^n \times G(m,n)} \phi(x,P) \chi_{{\rm spt\,}V}(x,P)\,dV(x,P)\\
& = & \int_{\R^n \times G(m,n)} \phi(x,P)\,dV(x,P).
\end{eqnarray*}
Thus we can replace the $\mu_P$ by $\mu_P \llcorner K_P$ without change.
Likewise, we may assume that $\mu_{P,z}$, for each $P \in G(m,n)$ and 
$z \in P^\perp$, is supported in
$K_{P,z} = {\rm spt\,}\mu_P \cap (z+P)$. In particular, for 
$x \in z+P$ we then have the implications
$$
x \in {\rm spt\,}\mu_{P,z}\, \Rightarrow\,  x \in {\rm spt\,}\mu_P 
\,\Rightarrow\, (x,P) \in {\rm spt\,}V.
$$}
\end{remark}

\begin{lemma} \label{lemmatest} Let $V$ be an $m$-varifold in $\R^n$ with 
support in $\overline{\Omega}$, which is orthogonal to $S$ along $\Gamma$
and has curvature $B_V = 0$. Then for $\nu$-a.e. $P \in G(m,n)$ and $\mu_P^\perp$-a.e.
$z \in P^\perp$ 
%the following holds: 
%up to null sets $E \subset G(m,n)$ and 
%$E_p \subset P^\perp$ for $\nu$ and $\mu_P^\perp$ respectively,
$$
\int_{z+P} \langle D\varphi(x),P \rangle\,d\mu_{P,z}(x) = 0
\quad \mbox{ for all $\varphi \in C^1(\R^n,\R^n)$ 
tangential on $S$.}
$$
\end{lemma}

{\em Proof. }We use the ansatz $\phi(x,P) = \varphi(x)\psi(P)$ 
in the orthogonal boundary condition, see Definition 
\ref{defadmissible}. This gives, since $\varphi(x)$
is tangential by assumption
\begin{eqnarray*}
0 & = & \int_{\R^n \times G(m,n)} \langle D_x \phi(x,P),P \rangle\,dV(x,P)\\
& = & \int_{G(m,n)} \psi(P) 
\int_{\R^n} \langle D\varphi(x),P \rangle\,d\mu_P(x)\,d\nu(P).
\end{eqnarray*}
Thus for any $P \in G(m,n)$, up to a set $E(\varphi)$ with measure
$\nu(E(\varphi)) =0$, we have
$$
\int_{\R^n} \langle D\varphi(x),P \rangle\,d\mu_P(x) = 0.
$$
We may also apply the above equation with $\varphi(x)$ replaced by $\sigma(P^\perp x) \varphi(x)$ (for a smooth function $\sigma$ and $\varphi$ exactly as before). We have
$$
D_\tau (\sigma(P^\perp x)) = D\sigma(P^\perp x) P^\perp\tau = 0 
\quad \mbox{ for }\tau \in P.
$$
Thus we get
$$
0 = \int_{P^\perp} \sigma(z) 
\int_{z+P} \langle D\varphi(x),P \rangle\,d\mu_{P,z}(x)\,d\mu_P^\perp(z).
$$
For any $z \in P^\perp$, up to a set $E_P(\varphi)$ with 
$\mu_P^\perp(E_P(\varphi)) = 0$, we conclude
$$
0 = \int_{z+P} \langle D\varphi(x),P \rangle\,d\mu_{P,z}(x).
$$
This is the desired result, except that the null sets depend 
on the test function $\varphi$. However, the set of $C^1$ functions 
$\varphi$  which are tangential along 
$\partial \Omega$ is separable. We can choose fixed null sets 
such that the identity holds for $\varphi$ in a countable
dense set, and hence for all $\varphi$ by approximation.
\kasten

\begin{lemma}[zero curvature] \label{lemmazero} 
Let $V$ be a varifold in $\Omega$ with curvature $B = 0$. Then 
for $\nu$-a.e. $P \in G(m,n)$ and $\mu_P^{\perp}$-a.e. 
$z \in P^\perp$, the following holds: for any connected 
component $U$ of $\Omega \cap (z+P)$ one has 
\begin{equation}
\label{eqzero} 
\mu_{P,z} = \theta_U\,{\cal L}^m \mbox{ on }U, \quad \mbox{ for a constant }
\theta_U \geq 0.
\end{equation}
\end{lemma}
The lemma follows from the well-known constancy lemma below,
using test functions $\chi(x) = \omega(x) \tau_\alpha$ 
for an orthonormal basis $\tau_\alpha$ of $P$.

\begin{lemma}[constancy] \label{lemmaconstancy}
Let $\mu$ be a Radon measure on an open domain $U \subset \R^m$, 
and assume that for all $\omega \in C^\infty_c(U)$ 
$$
\int_U \partial_\alpha \omega(x)\,d\mu(x) = 0 \quad \mbox{ for }
\alpha = 1,\ldots,m.
$$
Then $\mu = \theta\,{\cal L}^m$ on $U$ for some $\theta \in [0,\infty)$. 
\end{lemma}

The proof involving mollification is standard (cf. also Theorem 26.27 in \cite{Sim84}). 

%{\em Proof. }Let $\eta \in C^\infty_c(B_1(0))$ and 
%$\eta_\ve(x) = \ve^{-m} \eta(\frac{x}{\ve})$. For $\varphi \in C^0_c(U)$
%and sufficiently small $\ve > 0$, we compute
%$$
%\int_U \eta_\ve \ast \varphi(x)\,d\mu(x)
%= \int_U \int_{\R^m} \eta_\ve(x-y) \varphi(y)\,dy\,d\mu(x)
%= \int_{\R^m} \theta_\ve(y) \varphi(y)\,dy,
%$$
%where $\theta_\ve(y) = \int_U \eta_\ve(x-y)\,d\mu(x)$.
%Inserting $\varphi = \partial_\alpha \omega$, we obtain by assumption
%$$
%\int_{\R^n} \theta_\ve(y) \partial_\alpha \omega(y)\,dy 
%= \int_U \eta_\ve \ast \partial_\alpha \omega(x)\,d\mu(x)
%= \int_U \partial_\alpha (\eta_\ve \ast \omega)(x)\,d\mu(x)
%= 0.
%$$
%Thus $\theta_\ve \in [0,\infty)$ is constant in $U$, and   
%$$
%\int_U \eta_\ve \ast \varphi(x)\,d\mu(x) 
%= \theta_\ve \,\int_{\R^n} \varphi(y)\,dy.
%$$
%Chosing $\varphi$ with nonzero integral, we get that $\theta_\ve$
%converges to a constant $\theta \in [0,\infty)$ as $\ve \searrow 0$.
%This shows $\mu = \theta\,{\cal L}^m$ on $U$. \kasten

\begin{lemma}[dichotomy] \label{lemmadichotomy}
Let $V$ be an $m$-varifold in $\R^n$ with support in 
$\overline{\Omega}$, which is orthogonal to $S$ along
$\Gamma$ and has curvature $B_V = 0$. Then for any 
$(x_0,P_0) \in {\rm spt\,}V$ with $x_0 \in S$ 
one has the alternative 
\begin{equation}
\label {eqdichotomy}
\nu^S(x_0) \in P_0 \quad \mbox{ or } \quad P_0 \subset T_{x_0} S.
\end{equation}
In the second case, there exists $\varrho > 0$ such that 
\begin{equation}
\label{eqflatpiece}
(x,P_0) \in {\rm spt\,}V \quad \mbox{ for all }
x \in (x_0 + P_0) \cap B_\varrho(x_0). 
\end{equation}
In particular, $(x_0 + P_0) \cap B_\varrho(x_0) \subset \overline{\Omega}$.
\end{lemma}

{\em Proof. } Let $(x_0,P_0) \in {\rm spt\,}V$, $x_0 \in S$, and assume
\begin{equation}
\label{eqnotorthogonal}
\nu^S(x_0) \notin P_0.
\end{equation}
By continuity, there exists $r_0 > 0$ such that for all 
$0 < \delta,\varrho \leq r_0$ we have
$$
\nu^S(x) \notin P \quad \mbox{ for all } x \in S \cap B_{2\varrho}(x_0),\,
P \in B_\delta(P_0).
$$
It follows that 
$$
\nu^S(x) \wedge Q \notin B_\delta(P_0) \quad \mbox{ for all }
x \in S \cap B_{2\varrho}(x_0),\,Q \in G_{m-1}(T_xS).
$$
Thus for all $\varphi \in C^1_c(B_{2\varrho}(x_0),\R^n)$ and all 
$\psi \in C^1_c(B_\delta(P_0))$, we have 
\begin{eqnarray*}
&& \int_{S} \langle \nu^S(x),\varphi(x) \rangle \int_{G_{m-1}(T_x S)}
\psi(\nu^S(x) \wedge Q)\,d\Gamma^x(Q)\,d\sigma_{\Gamma}(x)\\
& = &  \int_{S \cap B_{2\varrho}(x_0)} \hspace{-5mm} \langle \nu^S(x),\varphi(x) 
\rangle \int_{G_{m-1}(T_x S)}
\underbrace{\psi(\nu^S(x) \wedge Q)}_{= 0}\,d\Gamma^x(Q)\,d\sigma_{\Gamma}(x) =0.
\end{eqnarray*}
Definition \ref{defadmissible}, with $B_V = 0$ and 
$\phi(x,P) = \varphi(x) \psi(P)$, now implies 
$$
0 = \int_{G(m,n)} \psi(P) \int_{\R^n} 
\langle P,D\varphi(x) \rangle\,d\mu_P(x)\,d\nu(P).
$$
Thus for $\nu$-a.e. $P \in B_\delta(P_0)$ we obtain
$$
\int_{\R^n} \langle P,D\varphi(x) \rangle\,d\mu_P(x) = 0
\quad \mbox{ for all }\varphi \in C^1_c(B_{2\varrho}(x_0),\R^n).
$$
Fix such a $P \in B_\delta(P_0)$, and consider a test 
function $\varphi(x) = \chi(Px) \zeta(P^\perp x)$, 
where $\chi \in C^1_c(P \cap B_{\varrho}(Px_0),\R^n)$ and 
$\zeta \in C^1_c(P^\perp \cap B_{\varrho}(P^\perp x_0))$. 
Then 
$$
\int_{P^\perp} \zeta(z) \int_{z+P} 
\mathrm{div}_P \chi(Px)\,d\mu_{P,z}(x) \,d\mu_P^\perp(z) = 0.
$$
Therefore we conclude, for $\mu_P^\perp$-a.e. 
$z \in P^\perp \cap B_{\varrho}(P^\perp x_0)$, 
\begin{eqnarray*}
0 & = & \int_{z+P} \mathrm{div}_P \chi(Px) \,d\mu_{P,z}(x)\\ 
& = & \int_{z+P} \mathrm{div}_P \chi(x-z) \,d\mu_{P,z}(x)
\quad \mbox{ for all }\chi \in C^1_c(P \cap B_{\varrho}(Px_0),\R^n).
\end{eqnarray*}
Thus for $\nu$-a.e. $P \in B_\delta(P_0)$ and $\mu_P^\perp$-a.e. 
$z \in P^\perp \cap B_{\varrho}(P^\perp x_0)$, we have 
\begin{equation}
\label{eqlocalconstancy}
\mu_{P,z} = \theta\,{\cal L}^m \llcorner (z+P) \quad
\mbox{ on }\,z+P \cap B_{\varrho}(Px_0) \quad 
\mbox{ for some }\theta \in [0,\infty).
\end{equation} 
%Now assume by contradiction that $P_0 \backslash T_{x_0}S \neq \emptyset$.
%Then $x_0$ is a regular point of the distance function $d_S|_{P_0}$, in fact
%$$
%\nabla^{P_0} d_S(x_0) = P_0 \nabla d_S(x_0) = 
%- P_0 \nu^S(x_0) \neq 0.
%$$
%Therefore $(x_0+P_0) \cap \Omega$ is a domain of class
%$C^1$ near its boundary point $x_0$. In particular, 
%$(x_0 + P_0) \cap B_\ve(x_0)$ has nonempty intersection .
%with $\R^n \backslash \overline{\Omega}$. 
%But then for $\ve > 0$ sufficiently small, there is also 
%a nonempty intersection with $\R^n \backslash \overline{\Omega}$
%of the set 
%$$
%z +  (P \cap B_\ve(y_0)) = z-z_0 + (x_0 + P) \cap B_\ve(x_0),
%$$
%for any $P \in B_\ve(P_0)$, $y_0 = Px_0$, $z_0 = P^\perp x_0$
%and $z \in P^\perp \cap B_\ve(z_0)$. This just follows by 
%continuity. We now conclude from (\ref{eqlocalconstancy}) 
%that $\mu_{P,z}$ is zero on $z-z_0 + (x_0+P) \cap B_{\ve/2}(x_0)$, 
%for all $P$, $z$ as above. This means that $(x_0,P_0)$ 
%is not in the support of the varifold $V$, contradicitng 
%the assumption.\\
%\\
We claim that for any $y_0 \in P_0 \cap B_{\varrho/2}(0)$, we 
have $(x_0+y_0,P_0) \in {\rm spt\,}V$. Let $B_\sigma(x_0+y_0)$
be given, wlog with $\sigma \leq \frac{\varrho}{4}$.
For $x \in B_\sigma(x_0+y_0)$ we have  
$$
|Px-Px_0| \leq |P(x-(x_0+y_0))| + |(P-P_0)y_0| + |y_0|
\leq \frac{3\varrho}{4} + \frac{\delta \varrho}{2},
$$
$$
|P^\perp x-P^\perp x_0| \leq
|P^\perp(x-(x_0+y_0))| + |(P^\perp-P^\perp_0)y_0|
\leq \frac{\varrho}{4} + \frac{\delta \varrho}{2}.
$$
Thus $|P(x-x_0)|,\,|P^\perp(x-x_0)| < \varrho$ for all 
$x \in B_\sigma(x_0+y_0)$, provided $\delta < \frac{1}{4}$.
Now we consider a ball $B_\ve(x_0+Py_0)$ and estimate
$$
|(x_0+y_0)-(x_0+Py_0)| = |P^\perp y_0| = |(P^\perp-P_0^\perp)y_0| 
\leq \frac{\delta \varrho}{2}.
$$
We choose $\ve = \frac{\sigma}{2}$, and conclude
$$
B_\ve(x_0+Py_0) \subset B_\sigma(x_0+y_0) \quad 
\mbox{ for any } \delta < \frac{\sigma}{\varrho}.
$$
Now we compute for fixed $P \in B_\delta(P_0)$ satisfying \eqref{eqlocalconstancy}
\begin{eqnarray*}
\mu_P(B_\sigma(x_0+y_0)) & \geq &  \mu_P(B_\ve(x_0+Py_0))\\
& = & \int_{P^\perp} \mu_{P,z}\big((z+P) \cap B_\ve(x_0+Py_0)\big) \,d\mu_P^\perp(z)\\
& = & \int_{P^\perp} \theta_{P,z} {\cal L}^m\big((z+P) \cap B_\ve(x_0+Py_0)\big)\,d\mu_P^\perp(z)\\
& = & \int_{P^\perp} \theta_{P,z} {\cal L}^m\big((z+P) \cap B_\ve(x_0)\big)\,d\mu_P^\perp(z)\\
& = & \mu_P\big(B_\ve(x_0)\big).
\end{eqnarray*}
Integrating over all $P \in B_\delta(P_0)$ we conclude
\begin{eqnarray*}
V\big(B_\sigma(x_0+y_0) \times B_\delta(P_0)\big) 
& = & \int_{B_\delta(P_0)} \mu_P(B_\sigma(x_0+y_0))\,d\nu(P)\\
& \geq & \int_{B_\delta(P_0)} \mu_P(B_\ve(x_0))\,d\nu(P)\\
& = & V\big(B_\ve(x_0) \times B_\delta(P_0)\big) > 0,
\end{eqnarray*}
since $(x_0,P_0) \in {\rm spt\,}V$.
This proves $(x_0 + y_0,P_0) \in {\rm spt\,}V$ for any 
$y_0 \in P_0 \cap B_{\varrho/2}(0)$. As a direct consequence, 
it follows that $(x_0 + P_0) \cap B_{\varrho/2}(x_0) \subset 
\overline{\Omega}$. \kasten 

\begin{lemma} \label{lemmaboundaryconstancy} 
Let $V$ be an $m$-varifold in $\R^n$ with support
in $\overline{\Omega}$, which is orthogonal to $S$ along $\Gamma$ 
and has curvature $B_V = 0$. Let $x_0 \in S$, and assume 
$(x_0,P_0) \in {\rm spt\,}V$ where $P_0 \subset T_{x_0}S$. 
Then for $z_0 = P_0^\perp x_0$ we have near $x_0$, if $P_0,z_0$ is as in Lemma \ref{lemmatest}, 
$$
\mu_{P_0,z_0} = \theta \,{\cal L}^m \llcorner (z_0+P_0) 
\quad \mbox{ for some constant }\theta \geq 0.
$$
\end{lemma}

{\em Proof. } For any $\varphi \in C^1(\R^n,\R^n)$ which is tangential
on $S$, we have by Lemma \ref{lemmatest} 
%here exist null sets $E \subset G(m,n)$ and $E_P \subset P^\perp$
%with respect to $\nu$ and $\mu_P^\perp$ respectively, such that unless 
%$P_0 \in E$ or $z_0 \in E_{P_0}$ for any 
\begin{equation}
\label{eqtest}
\int_{z_0+P_0} \langle D\varphi(x),P_0 \rangle \,d\mu_{P_0,z_0}(x) = 0.
\end{equation}
Choose a neighborhood $U \times I$ of
$x_0$ on which $\Omega$ has a subgraph representation
$$
\Omega \cap (U \times I) = \{(y,z) \in U \times I: z < u(y)\}.
$$
After rotation we assume $x_0 = 0$, $T_{x_0}S = \R^{n-1}$ 
and $P_0 = \R^m \times \{0\}$ where $m \leq n-1$. For any 
$\sigma \in C^1_c(U \times I)$ we consider the test function
$$
\varphi(y,z) = \sigma(y,z)\,\big(e_\alpha + \partial_\alpha u(y)e_n\big)
\quad \mbox{ where }x = (y,z) \mbox{ and }1 \leq \alpha \leq m. 
$$
By construction the vector field $\varphi$ is tangential along $\partial \Omega$. 
%$\big\langle e_\alpha + \partial_\alpha u(y)e_n, (-\nabla u(y),1) \big\rangle = 0$.
Furthermore
$$
\langle D\varphi(x),P_0 \rangle = 
\sum_{\beta=1}^m \langle \partial_\beta \varphi(x), e_\beta \rangle
= \partial_\alpha \sigma(x).
$$
Therefore (\ref{eqtest}) implies 
$$
\int_{\R^m} \partial_\alpha \sigma(x)\,d\mu_{\R^m,0} =0 \quad 
\mbox{ where }(P_0,z_0) = (\R^m,0).
$$
Since any function with compact support in $\R^m \cap U$ can be 
extended to a function with compact support in $U \times I$, the 
claim follows from Lemma \ref{lemmaconstancy}. \kasten

\begin{theorem}[Orthogonal Varifolds with Curvature Zero] 
\label{thmorthozero} Let $V$ be an $m$-varifold in
$\R^n$ with support in $\overline{\Omega}$, which is orthogonal
to $S = \partial \Omega$ along an $(m-1)$-varifold $\Gamma$ and has 
curvature $B_V = 0$. Then for all $\phi \in C^0(\R^n \times G(m,n))$ 
and $\varphi \in C^0(G_{m-1}(TS))$, respectively, the following 
representations hold:
\begin{eqnarray*}
%\label{eqorthozero1}
V(\phi) & = & \int_{G(m,n)} \int_{P^\perp} \sum_{\Delta \in S^m_\perp(\Omega,z+P)}
\theta_\Delta \int_\Delta \phi(x,P)\,d{\cal L}^m(x)\,
d\mu_P^\perp(z)\,d\nu(P),\\
%\label{eqorthozero2}
\Gamma(\varphi) & = & \int_{G(m,n)} \int_{P^\perp} 
\sum_{\Delta \in S^m_\perp(\Omega,z+P)}
\theta_\Delta \int_{\partial \Delta}
\varphi(x,P \cap T_x S)\,d{\cal H}^{m-1}(x)\,d\mu_P^\perp(z)\,d\nu(P).
\end{eqnarray*}
Here $S^m_\perp(\Omega,z+P)$ is the finite, possibly empty 
set of orthogonal $m$-slices of $\Omega$ in the affine plane $z+P$,
and the $\theta_\Delta \geq 0$ are constants.
\end{theorem}

{\em Proof. }Let $\Delta$ be any $m$-slice of $\Omega$ in an affine 
plane $z+P$, where $P \in G(m,n)$ and $z \in P^\perp$. By definition,
for any $x \in \Delta$ there exists $\varrho > 0$ such that
$$
(z+P) \cap B_\varrho(x) \subset \overline{\Omega}.
$$
In particular if $x \in \Delta \cap S$ then $P \subset T_{x}S$. Lemma
\ref{lemmazero} and Lemma \ref{lemmaboundaryconstancy} imply that
$\mu_{P,z} = \theta\, {\cal L}^m \llcorner (z+P)$ locally 
near any $x \in \Delta$. Since $\Delta$ is connected, we conclude
$$
\mu_{P,z} \llcorner \Delta = \theta_\Delta\, {\cal L}^m \llcorner \Delta 
\quad \mbox{ for some }\theta_\Delta \geq 0. 
$$
Actually, Lemma \ref{lemmaboundaryconstancy} requires that 
$(x,P) \in {\rm spt\,}V$, but otherwise the statement follows 
with $\theta_\Delta = 0$, see Remark \ref{rmksupports}. 
Now assume that $\theta_\Delta > 0$, and consider a boundary 
point $x \in \partial \Delta \subset S$. We then have 
$(x,P) \in {\rm spt\,}V$ again by Remark \ref{rmksupports}.
We claim that $\nu^S(x) \in P$. Otherwise Lemma 
\ref{lemmadichotomy} yields that
$(x+P) \cap B_\varrho(x) \subset \overline{\Omega}$ for 
some $\varrho > 0$. But $\Delta \cap B_\varrho(x)$ is 
nonempty, thus $(x+P) \cap B_\varrho(x) \subset \Delta$
and $x \in \Delta$, a contradiction. This 
shows that if $\theta_\Delta > 0$ then $\nu^S(x) \in P$ 
for all $x \in \partial \Delta$, in other words $\Delta$
is an orthogonal $m$-slice. Along $\partial \Delta$ 
we have
$$
\nabla ({\rm dist}_S|_{z+P})(x)
= P\, \nabla {\rm dist}_S(x) = \nu^S(x) \neq 0,
$$
therefore $\Delta$ is a bounded domain in $z+P$ of class $C^2$.
We denote by $\Delta_{P,z}$ the union of the finitely many 
orthogonal slices $\Delta \subset z+P$ with $\theta_\Delta > 0$,
and introduce the function
$$
\theta_{P,z}:z+P \to [0,\infty),\,\theta_{P,z}(x) = \theta_\Delta 
\mbox{ for }x \in \overline{\Delta}.
$$
%We now classify the curvature zero $m$-varifolds 
%in $\overline{\Omega}$ with orthogonal boundary condition. 
%For any $P \in G(m,n)$ and $z \in P^\perp$, let 
%$C_{P,z}$ be the set of points $x \in z+P$ such that 
%$\mu_{P,z}$ is not a constant multiple of 
%${\cal L}^m$ on $z+P$ near $x$. By definition
%$C_{P,z}$ is compact. We have shown that
%$$
%C_{P,z} \subset S,\, \mbox{ moreover }\nu_S(x) \in P 
%\mbox{ for any }x \in C_{P,z}.
%$$
%Thus at any $x \in C_{p,z}$ we compute
%$$
%\nabla ({\rm dist}_S|_{z+P})(x)
%= P(x) \nabla {\rm dist}_S(x) = \nu_S(x) \neq 0.
%$$
%This means that $(z+P) \cap \Omega$ is locally a 
%$C^1$ domain near $x$, with boundary $(z+P) \cap S$. 
%For an appropriate neighborhood $U \subset z+P$ of 
%$x$ and some constant $\theta > 0$, we get 
%$$
%\mu_{P,z} = \begin{cases}
%\theta\, {\cal L}^m & \mbox{ on }U \cap \Omega,\\
%0             & \mbox{ on }U \backslash \overline{\Omega}.
%\end{cases}
%$$
%In particular $U \cap S$ is also in $C_{P,z}$, and 
%$C_{P,z}$ is an open subset of $(z+P) \cap S$. In summary, 
%we conclude that $C_{P,z}$ is a compact, $(m-1)$-dimensional
%submanifold of $z+P$, possibly with several components, and 
%that $z+P$ and $S$ intersect orthogonally along $C_{P,z}$. 
The measure $\tilde{\mu}_{P,z} = {\cal L}^m \llcorner \theta_{P,z}$
satisfies, for any $\varphi \in C^1(z+P,P)$, 
\begin{eqnarray*}
\int_{z+P} \langle D\varphi(x), P \rangle\,d\tilde{\mu}_{P,z}(x) & = & 
\sum_{\Delta \in S^m_\perp(\Omega,z+P)} \theta_\Delta 
\int_\Delta \langle D\varphi(x), P \rangle\,d{\cal L}^m(x)\\
& = & - \sum_{\Delta \in S^m_\perp(\Omega,z+P)} \theta_\Delta
\int_{\partial \Delta} \langle \varphi(x),\nu^S(x) \rangle\,d{\cal H}^{m-1}(x).
%& = & -\int_{\partial \Delta_{P,z}} 
%\theta_{P,z}(x) \langle \varphi(x),\nu_S(x) \rangle \,d{\cal H}^{m-1}(x).
\end{eqnarray*} 
In other words we have, in the sense of distributions on $z+P$, 
\begin{equation}
\label{eqboundarymeasure}
\nabla \tilde{\mu}_{P,z} = 
- \theta_{P,z} ({\cal H}^{m-1} \llcorner \partial \Delta_{P,z}) \otimes \nu^S.
\end{equation}
%We finally remark that $C_{P,z}$ is empty if and only if $\mu_{P,z} =0$. 
%Namely if $\mu_{P,z}$ is a constant multiple of ${\cal L}^m \llcorner (z+P)$ 
%for any $x \in z+P$, then $\mu_{P,z} = \theta \,{\cal L}^m \llcorner (z+P)$
%globally for a constant $\theta \geq 0$. As $\mu_{P,z}$ has compact 
%support, we obtain $\theta = 0$. Thus for $V \neq 0$, we always have 
%some orthogonal intersection $S \cap (z+P)$ along $C_{P,z}$, for some 
%$P \in G(m,n)$ and $z \in P^\perp$. This statement is not affected 
%by the exceptional null sets, since they do not contribute to $V$.\\
On functions $\varphi \in C^0_c(\R^n \times G(m,n))$, we now consider the 
$m$-varifold $\tilde{V}$ given by 
\begin{equation}
\label{eqtildeV}
\tilde{V}(\varphi) = 
\int_{G(m,n)} \int_{P^\perp} \int_{z+P} 
\varphi(x,P)\,d\tilde{\mu}_{P,z}(x)\,d\mu_P^\perp(z)\,d\nu(P). 
\end{equation}
For $\varphi \geq 0$ we can estimate
\begin{eqnarray*}
0 \leq \tilde{V}(\varphi) & = &
\int_{G(m,n)} \int_{P^\perp} \int_{\Delta_{P,z}}
\varphi(x,P)\,\,d\mu_{P,z}(x)\,d\mu_P^\perp(z)\,d\nu(P)\\
& \leq & \int_{G(m,n)} \int_{P^\perp} \int_{z+P}
\varphi(x,P)\, d\mu_{P,z}(x)\,d\mu_P^\perp(z)\,d\nu(P)\\
& = & V(\varphi).
\end{eqnarray*}
In particular $\tilde{V}$ is continuous, in fact we have 
$$
|\tilde{V}(\varphi)| \leq \tilde{V}(|\varphi|) \leq V(|\varphi|) 
\leq {\bf M}(V)\,\|\varphi\|_{C^0(\R^n \times G(m,n))}.
$$
It follows that $\Lambda(\varphi):= V(\varphi) - \tilde{V}(\varphi)$ 
is a Radon measure on $\R^n \times G(m,n)$. 
Moreover $\tilde{V}(\varphi) = V(\varphi)$ whenever $\varphi$ has 
compact support in $\Omega \times G(m,n)$. Therefore $\Lambda$ 
has support in $S \times G(m,n)$. Using (\ref{eqboundarymeasure}) 
we further have for $\phi \in C^1(\R^n \times G(m,n),\R^n)$
\begin{eqnarray*} 
&& \int_{\R^n \times G(m,n)} \langle D_x\phi(x,P),P \rangle\,d\tilde{V}(x,P)\\
& = & \int_{G(m,n)} \int_{P^\perp} \int_{z+P} 
\langle D_x \phi(x,P), P\rangle \,d\tilde{\mu}_{P,z}(x)\,d\mu_P^\perp(z)\,d\nu(P)\\
& = & - \int_{G(m,n)} \int_{P^\perp} 
\int_{\partial \Delta_{P,z}} \theta_{P,z}(x) 
\langle \nu^S(x),\phi(x,P) \rangle\,d{\cal H}^{m-1}(x)\,
d\mu_P^\perp(z)\,d\nu(P). 
\end{eqnarray*}
On the other hand, $V$ satisfies the orthogonal boundary constraint
$$
\int_{\R^n \times G(m,n)} \langle D_x\phi(x,P),P \rangle\,dV(x,P)
= - \int_{G_{m-1}(TS)} \langle \nu^S(x),\phi(x,\nu^S(x) \wedge Q) \rangle\,\,
d\Gamma(x,Q).
$$
We take a test function of the form $\phi(x,P) = \varphi(x)$ where $\varphi|_{S} = 0$.
Then 
\begin{eqnarray*}
0 & = & \int_{\R^n \times G(m,n)} \langle D\varphi(x),P \rangle\,dV(x,P) 
- \int_{\R^n \times G(m,n)} \langle D\varphi(x),P \rangle\,d\tilde{V}(x,P)\\
& = & \int_{S \times G(m,n)} \langle D\varphi(x),P \rangle\,d\Lambda(x,P)\\
& = & \int_{S \times G(m,n)} \langle D\varphi(x)\nu^S(x),\nu^S(x) \rangle\,d\Lambda(x,P)\\
& = & \int_S \Lambda^x(G(m,n))\,
\langle D\varphi(x)\nu^S(x),\nu^S(x) \rangle\,d\mu_\Lambda(x),
\end{eqnarray*}
where $\mu_\Lambda, \Lambda^x$ are the measures arising from the disintegration of $\Lambda$, cf. Theorem \ref{lemmadisintegration}.
Taking now $\varphi(x) = \mathrm{dist}_S(x)\,\nu^S(x)$ yields $D\varphi(x) \nu^S(x) = \nu^S(x)$
along $S$. We conclude
$$
0 = \int_S \Lambda^x(G(m,n))\,d\mu_\Lambda(x) = \Lambda (\R^n \times G(m,n)).
$$
This shows $\Lambda = 0$ and hence $V = \tilde{V}$, which is the desired 
representaton for $V$. Finally, we claim that the boundary varifold
$\Gamma$ in Definition \ref{defadmissible} has the representation 
$$
\Gamma(\varphi) =  \int_{G(m,n)} \int_{P^\perp}
\int_{\partial \Delta_{P,z}} \theta_{P,z}(x)\, \varphi(x,P \cap T_x S)\,d{\cal H}^{m-1}(x)\,
d\mu_P^\perp(z)\,d\nu(P).
$$
First note that $P \cap T_x S \in G_{m-1}(T_x S)$ as $\nu^S(x) \in P$ for 
$x \in \partial \Delta_{P,z}$. 
For given $\varphi \in C^1(G_{m-1}(TS))$ we consider as test 
function $\phi(x,P) = \varphi(x,P \cap T_x S) \nu^S(x)$, suitably 
extended to a $C^1$ function on $\R^n \times G(m,n)$. From $V = \tilde{V}$ 
and $B_V = 0$ we get
\begin{eqnarray*}
&& \int_{G_{m-1}(TS)} \varphi(x,Q)\,d\Gamma(x,Q)\\ 
& = & \int_{G_{m-1}(TS)} \langle \nu^S(x),\phi(x,\nu^S(x) \wedge Q) \rangle\,d\Gamma(x,Q)\\
& = & - \int_{\R^n \times G(m,n)} \langle D_x \phi(x,P),P \rangle\,dV(x,P)\\
& = & \int_{G(m,n)} \int_{P^\perp}
\int_{\partial \Delta_{P,z}} \theta_{P,z}(x) \langle \nu^S(x),\phi(x,P) \rangle\,d{\cal H}^{m-1}(x)\,
d\mu_P^\perp(z)\,d\nu(P)\\
& = &  \int_{G(m,n)} \int_{P^\perp} \int_{\partial \Delta_{P,z}} \theta_{P,z}(x)
\varphi(x,P \cap T_x S)\,d{\cal H}^{m-1}(x)\,d\mu_P^\perp(z)\,d\nu(P). 
\end{eqnarray*}\kasten

%For $\varphi(x,Q) =  \langle \nu_S(x),\phi(x, \nu_S(x)\wedge Q) \rangle$ 
%we compute, assuming $\nu_S(x) \in P$,
%$$
%\varphi(x,\nu_S(x) \llcorner P) = 
%\langle \nu_S(x),\phi\big(x, \nu_S(x)\wedge (\nu_S(x) \llcorner P)\big)\rangle
%= \langle \nu_S(x),\phi(x, P)\rangle.
%$$
%Recalling $\nu_S(x) \in P$ for $x \in C_{P,z}$ we obtain for this 
%choice of $\varphi(x,Q)$
%\begin{eqnarray*}
%\tilde{\Gamma}(\varphi) =  - \int_{G(m,n)} \int_{P^\perp}
%\int_{C_{P,z}} \theta_{P,z}(x) \langle \nu_S(x),\phi(x,P)\rangle\,
%\,d{\cal H}^{m-1}(x)\,d\mu_P^\perp(z)\,d\nu(P).
%\end{eqnarray*}
%It remains to show that $\tilde{\Gamma}(\varphi)$ defines a 
%an $(m-1)$-varifold on $S$.

\section{Mass and Topology Bounds} 

In this section Theorem \ref{thmorthozero} is applied to prove 
the bounds stated in the introduction. 

%{\bf Theorem 3 }Let $\Omega \subset \R^n$ be a bounded domain of 
%calls $C^2$. If $\Omega$ has no orthogonal $m$-slices, then
%$\varkappa^{m,p}(\Omega) > 0$ for any $p \in [1,\infty]$.
%Moreover for any varifold $V \in {\cal V}^p_m(\Omega)$ with 
%associated boundary varifold $\Gamma$, we have 
%\begin{eqnarray}
%\label{eqvarifoldmass} 
%{\bf M}(V) & \leq & \frac{1}{\varkappa^{m,p}(\Omega)}\,\|B_V\|^p_{L^p(V)},\\
%\label{eqboundarymass}
%{\bf M}(\Gamma) & \leq & C\,\|B_V\|^p_{L^p(V)} \quad \mbox{ where }
%C = C(n,m,p,\Omega).
%\end{eqnarray}

{\em Proof of Theorem \ref{thmvarifolds}}. We first prove that
$\varkappa^{m,p}(\Omega) > 0$. As 
$\varkappa^{m,p}(\Omega)^{\frac{1}{p}} \geq \varkappa^{m,1}(\Omega)$
by definition, it suffices to consider the case $p = 1$. 
Assume by contradiction that there is a sequence
$V_k \in {\bf CV}^m_{\perp}(\Omega)$ with
${\bf M}(V_k) = 1$, such that $\|B_{V_k}\|_{L^1(V_k)} \to 0$.
For the associated boundary varifolds $\Gamma_k$, we have 
by Lemma \ref{lemmamasses}
$$
{\bf M}(\Gamma_k) \leq C\,\big({\bf M}(V_k) + \|B_{V_k}\|_{L^1(V_k)}\big)
\leq C \quad \mbox{ where }C = C(n,\Omega).
$$
Thus we can assume that $V_k \to V$ where again ${\bf M}(V) = 1$, 
and also $\Gamma_k \to \Gamma$ as varifolds. By passing to 
the limit in (\ref{eqadmissible}) we get
$$
\int_{\R^n \times G(m,n)} \langle D_x \phi,P \rangle\,dV(x,P)
= - \int_{G_{m-1}(TS)}
\langle \nu^S(x),\phi(x,\nu^S(x) \wedge Q)\rangle\,d\Gamma(x,Q).
$$
This means that $V$ is orthogonal to $S$ along $\Gamma$ and
has curvature zero. But $\Omega$ has no orthogonal $m$-slices 
by assumption, and hence Theorem \ref{thmorthozero} implies
that $V = 0$, a contradiction. For given nonzero
$V \in {\bf CV}^m_\perp(\Omega)$, let $V' = \lambda V$ where
$\lambda = {\bf M}(V)^{-1}$, thus ${\bf M}(V') = 1$. By definition of 
$\varkappa^{m,p}(\Omega)$ and (\ref{eqscaling})
$$
\|B_V\|^p_{L^p(V)}  = {\bf M}(V)\,\|B_{V'}\|^p_{L^p(V')}
\geq {\bf M}(V)\, \varkappa^{m,p}(\Omega).
$$
The estimate for $\Gamma$ follows from Lemma \ref{lemmamasses} 
and Remark \ref{rmkyoung}. \kasten

%Theorem \ref{thm1} is just a special case of Theorem \ref{thmvarifolds},
%except for estimate (\ref{eqeuler}) which will be proved below
%when discussing the two-dimensional case.

{\em Proof of Theorem \ref{thm1}.} As discussed in Section 2 any $\Sigma \in {\bf S}^m_\perp(\Omega)$
with $\int_\Sigma |A_\Sigma|^p \, d\mu_\Sigma < \infty$ induces a curvature varifold in 
${\bf CV}^m_\perp(\Omega)$. Hence the mass bound (\ref{eqthm1bound}) follows 
from Theorem \ref{thmvarifolds}. Now we turn to the special case when $m = 2$ and $p = 2$.
Let $\Sigma \in S_\perp^2(\Omega)$ with $\int_\Sigma |A_\Sigma|^2 \, d\mu_\Sigma < \infty$. 
%Consider for $\varepsilon> 0$ the domains $\Omega_\varepsilon = \Omega \setminus \overline{U_\varepsilon^+(S)}$
%and put $S_\varepsilon = \partial \Omega_\varepsilon$. One readily checks that the 
%surface $\Sigma_\varepsilon = \Sigma \backslash U_\ve^{+}(S)$ is for $\ve > 0$ small 
%an open subset of class $C^2$ of ${\rm int}(\Sigma)$, and $\partial \Sigma_\varepsilon$ converges 
%to $\partial \Sigma$ in the $C^1$ topology. By \cite{AK16}, see the calculation 
%before (2.10) there, the geodesic curvature $\varkappa_{\partial \Sigma_\ve}$ of $\partial \Sigma_\ve$
%is given by
%\begin{equation}
%\label{eqgeodesiccurvature}
%\varkappa_{\partial \Sigma_\ve} = 
%\frac{h^{S_\ve}(\tau_\ve,\tau_\ve)}{\langle \eta_{\partial \Sigma_\ve}, \nu^{S_\ve} \rangle},
%\end{equation}
%where $\tau_\ve$ is a unit tangent of $\partial \Sigma_\varepsilon$. 
%The Gau{\ss}-Bonnet theorem applied to $\Sigma_\varepsilon$ yields 
%\begin{equation*}
%2\pi \chi(\Sigma_\varepsilon) = \int_{\Sigma_\varepsilon} K_\Sigma \, d \mu_\Sigma 
%+ \int_{\partial \Sigma_\varepsilon} 
%\frac{h^{S_\varepsilon}(\tau_\ve,\tau_\ve)}{\langle \eta_{\partial \Sigma_\ve},\nu^{S_\varepsilon}\rangle}
 % \, d\sigma_{\partial \Sigma_\ve}.
%\end{equation*}
%The boundary integral depends only on $C^1$ data of $\Sigma_\ve$, and 
%on $C^2$ data of $S_\ve$, moreover one has $\chi(\Sigma_\varepsilon) = \chi(\Sigma)$
%for $\ve > 0$ small. Using that $\eta_{\partial \Sigma} = \nu^S|_{\partial \Sigma}$ 
%and $|K_\Sigma| \leq \frac{1}{2} |A_\Sigma|^2 \in L^1(\mu_\Sigma)$,
%we can pass to the limit $\ve \searrow 0$ to obtain
Using $\eta_{\partial \Sigma} = \nu^S\vert_{\partial \Sigma}$ we have by the Gauss-Bonnet Theorem and Meusnier's Theorem 
\begin{equation}
\label{eqgaussbonnet}
2\pi \chi(\Sigma) = \int_{\Sigma} K_\Sigma \, d\mu_\Sigma 
+ \int_{\partial \Sigma} h^S(\tau,\tau)\, d\sigma_{\partial \Sigma}.
\end{equation}
This is immediate if $\Sigma$ is of class $C^2$ up to the boundary. For general $\Sigma \in S_\perp^2(\Omega)$ one uses an approximation, see also the calculation before (2.10) in \cite{AK16}.
Now \eqref{eqgaussbonnet} yields 
%Using further that $\Sigma_\varepsilon \subset \Sigma,$ $|\partial \Sigma_\varepsilon| \rightarrow |\partial \Sigma|$ ($\varepsilon \rightarrow 0$) and $|K_\Sigma| \leq  \frac{1}{2}|A_\Sigma|^2$ we obtain
$$
2\pi\, |\chi(\Sigma)|
\leq   \frac{1}{2} \int_{\Sigma} |A_\Sigma|^2\,d\mu_\Sigma
+ \|h^S\|_{C^0(S)}\, |\partial \Sigma|.
$$
The claimed estimate for $\chi(\Sigma)$ follows easily from this and equation (\ref{eqthm1bound}).
\kasten

{\em Proof of Theorem \ref{thmboundsnearboundary}.} To 
prove claim (\ref{eqboundsnearboundary}) it is sufficient to 
bound $|\Sigma|$, since then $|\partial \Sigma|$ is estimated
by Lemma \ref{lemmamasses}. Assume by contradiction 
that there is a sequence $\Sigma_k \in {\bf S}^m_\perp(\Omega)$
with $\Sigma_k \subset U_\delta(S)$, such that 
$$
\int_{\Sigma_k} |A_{\Sigma_k}|^p\,d\mu_{\Sigma_k}
< \frac{1}{k}\, |\Sigma_k|.
$$
Passing to $V_k = \Sigma_k/|\Sigma_k| \in {\bf CV}^m_\perp(\Omega)$ 
we have ${\bf M}(V_k) = 1$, and (\ref{eqscaling}) yields
$$
\|B_{V_k}\|_{L^p(V_k)}^p =
2^{p/2} \|A_{\Sigma_k}\|_{L^p(\Sigma_k)}^p/|\Sigma_k| 
< 2^{p/2} \frac{1}{k} \to 0.
$$
$V_k$ satisfies Definition \ref{defadmissible} with 
boundary $\Gamma_k = \partial \Sigma_k/|\Sigma_k|$. As 
${\bf M}(\Gamma_k) \leq C$ by Lemma \ref{lemmamasses}, we 
may assume that $V_k \to V$ and $\Gamma_k\to \Gamma$ as 
varifolds. It follows that ${\bf M}(V) = 1$, $V$ is orthogonal
to $S$ and has curvature zero. Moreover
\begin{equation}\label{muVspt}
{\rm spt\,}\mu_V \subset \overline{U_\delta(S)}. 
\end{equation}
Let $\Delta$ be any orthogonal $m$-slice of $\Omega$. Choose 
$p \in \Delta$ such that $\varrho = \mathrm{dist}(p,\partial \Delta)$
is maximal. Then there exist at least two points $x \in \partial \Delta$
such that $p = x + \varrho \nu^S(x)$, and thus the normal injectivity 
radius of $S$ satisfies $\varrho_S \leq \varrho$. We 
get for $\Omega_\delta = \Omega \backslash \overline{U_\delta(S)}$
$$
|\Delta \cap \Omega_\delta| \geq c(m) (\varrho_S-\delta)^m > 0.
$$
On the other hand, by comparing with a ball containing $\Omega$, 
we have $|\Delta| \leq C$ for $C = C(m,\Omega)$. From the 
representation in Theorem \ref{thmorthozero} we now get
\begin{eqnarray*}
\mu_V(\Omega_\delta) & = & 
\int_{G(m,n)} \int_{P^\perp} \sum_{\Delta \in S^m_\perp(\Omega,z+P)}
\theta_\Delta\, |\Delta \cap \Omega_\delta|\,d\mu_P^\perp(z)\,d\nu(P)\\
& \geq & c(m,\Omega) (\varrho_S-\delta)^m 
\int_{G(m,n)} \int_{P^\perp} \sum_{\Delta \in S^m_\perp(\Omega,z+P)}
\theta_\Delta\, |\Delta|\, d\mu_P^\perp(z)\,d\nu(P)\\
& = & c(m,\Omega) (\varrho_S-\delta)^m\, M(V) > 0.
\end{eqnarray*}
This contradicts \eqref{muVspt} and proves claim (\ref{eqboundsnearboundary}). 
The bound on the Euler characteristic in the case $m = p = 2$ follows 
as in the proof of Theorem \ref{thm1}.
\kasten

\begin{example} {\em If $\Omega$ is a round ball or an $n$-dimensional 
ellipsoid, then any orthogonal slice $\Delta$ contains the origin. 
Adapting the arguments yields the following: for any $\ve > 0$ 
there is a constant $C = C(n,m,p,\Omega,\ve) < \infty$ such that 
for any $\Sigma \in S^m_\perp(\Omega)$ 
$$
|\Sigma| + |\partial \Sigma| \leq C \int_\Sigma |A_\Sigma|^p\,d\mu_\Sigma
\quad \mbox{ whenever }\Sigma \subset \Omega \backslash B_\ve(0).
$$}
\end{example}

%& \leq & C(n,\Omega)\,\int_{\Sigma} |A_\Sigma|^2\,d\mu_\Sigma.
%We remark further that if $h^S \geq k > 0$, then again by Gau{\ss}-Bonnet
%\begin{eqnarray*}
%k\,|\partial \Sigma| + 2\pi |\chi(\Sigma)^{-}| 
%& \leq & \int_{\partial \Sigma} \varkappa_g\,ds 
%+ 2\pi \chi(\Sigma)^{+} - 2\pi \chi(\Sigma)\\
%& = & 2\pi \chi(\Sigma)^{+} - \int_{\Sigma} K_{\Sigma}\,d\mu_\Sigma\\
%& \leq & 2\pi \chi(\Sigma)^{+} 
%+ \frac{1}{2} \int_\Sigma |A_\Sigma|^2\,d\mu_\Sigma.
%\end{eqnarray*}
%The present paper was motivated by the question whether the 
%condition of strict convexity could be relaxed.

%{\bf Theorem 2 } 
%Let $\Omega \subset \R^n$
%be a bounded domain of class $C^2$. Assume that whenever $P$ is
%an orthogonal $2$-plane for $\Omega$, the set
%$\overline{\Omega} \cap P$ is a topological disk.
%Then there is a constant $C < \infty$ depending on $n$
%and $\Omega$ such that
%\begin{equation}
%|\Sigma| + |\partial \Sigma| + |\chi(\Sigma)| \leq
%C\,\Big(\int_\Sigma |A_\Sigma|^2\,d\mu_\Sigma + \chi(\Sigma)^{+}\Big) \quad
%\mbox{ for any }\Sigma \in {\bf S}_2(\Omega).
%\end{equation}

{\em Proof of Theorem \ref{thmdisktype}. }It suffices to prove the bound for 
$|\Sigma|$. Assume by contradiction that there are immersed
surfaces $\Sigma_k \in {\bf S}^2_\perp(\Omega)$ such that
%since the bounds for $|\partial \Sigma|$ and $|\chi(\Sigma)|$ 
%follow then by Lemma \ref{lemmamasses} and by the argument in Corollary 
%\ref{corollarytopology}. Assume by contradiction that there
$$
\int_{\Sigma_k} |A_{\Sigma_k}|^2\,d\mu_{\Sigma_k} + \chi(\Sigma_k)^{+}
< \frac{1}{k}\, |\Sigma_k|.
$$
The sequence $V_k = \Sigma_k/|\Sigma_k| \in {\bf CV}^2_\perp(\Omega)$ 
has norm ${\bf M}(V_k) = 1$ and boundary $\Gamma_k = \partial \Sigma_k/|\Sigma_k|$.
After passing to a subsequence we have $V_k \to V$ and $\Gamma_k\to \Gamma$
as varifolds, where Definition \ref{defadmissible} holds in the 
limit with curvature $B_V = 0$. Now by (\ref{eqgaussbonnet}) we have,
denoting by $\tau_{\partial \Sigma_k}$ the unit 
tangent vector of $\partial \Sigma_k$,
$$
\int_{\partial \Sigma_k} 
h^S\big(\tau_{\partial \Sigma_k},\tau_{\partial \Sigma_k}\big)
\,d\sigma_{\partial \Sigma_k}
\leq 
2\pi \chi(\Sigma_k)^{+} - \int_{\Sigma_k} K_{\Sigma_k}\,d\mu_{\Sigma_k}.
$$
We get a well-defined function $\varphi:G_1(TS) \to \R$ by putting
$\varphi(x,Q) = h^S(x)(\tau,\tau)$ where $\tau \in Q$ with $|\tau| = 1$.
By the varifold convergence $\Gamma_k \to \Gamma$, we have
\begin{eqnarray*}
\int_{G_1(TS)} \varphi(x,Q)\,d\Gamma(x,Q) & = & 
\lim_{k \to \infty} \int_{G_1(TS)} \varphi(x,Q)\,d\Gamma_k(x,Q)\\
& = & \lim_{k \to \infty} \frac{1}{|\Sigma_k|}
\int_{\partial \Sigma_k}
h^S\big(\tau_{\partial \Sigma_k},\tau_{\partial \Sigma_k}\big)
\,d\sigma_{\partial \Sigma_k}\\
& \leq & \lim_{k \to \infty} \frac{1}{|\Sigma_k|}
\Big(2\pi \chi(\Sigma_k)^{+} 
- \int_{\Sigma_k} K_{\Sigma_k}\,d\mu_{\Sigma_k}\Big)\\
& = & 0.
\end{eqnarray*}
On the other hand Theorem \ref{thmorthozero} yields the representation
$$
\Gamma(\varphi) =  \int_{G(2,n)} \int_{P^\perp} \sum_{\Delta \in S^2_\perp(\Omega,z+P)}
\theta_\Delta \int_{\partial \Delta} \varphi(x,P \cap T_x S)\,d{\cal H}^1(x)\,
d\mu_P^\perp(z)\,d\nu(P).
$$
Since $\Delta$ is an orthogonal slice, we have by the computation above
$$
\varphi(x,P \cap T_x S) = 
h^S(x)\big(\tau_{\partial \Delta}(x), \tau_{\partial \Delta}(x)\big)
= \varkappa_{\partial \Delta}(x). 
$$
Here $\Delta$ is a domain in the plane $P$ and $\varkappa_{\partial \Delta}$ 
is the Euclidean curvature of its boundary. Now by assumption each
$\Delta$ is a topological disk, therefore by the Hopf Umlaufsatz
\cite{Hop89} 
$$
\int_{\partial \Delta} \varphi(x,P \cap T_x S)\,d{\cal H}^1(x)
= \int_{\partial \Delta} \varkappa_{\partial \Delta}(x)\,d{\cal H}^1(x)
= 2\pi.
$$
Inserting we obtain
$$
\Gamma(\varphi) = 2\pi \int_{G(2,n)} \int_{P^\perp} 
\sum_{\Delta \in S^2_\perp(\Omega,z+P)} \theta_\Delta\,
d\mu_p^\perp(z) \,d\nu(P).
$$
But this yields, since $|\Delta| \leq C(n,\Omega)$,
$$
1 = {\bf M}(V) = \int_{G(2,n)} \int_{P^\perp} 
\sum_{\Delta \in S^2_\perp(\Omega,z+P)} \theta_\Delta\,|\Delta|
\leq C(n,\Omega)\, \Gamma(\varphi) \leq 0.
$$
\kasten 

\begin{example} {\em The assumptions in Theorem
\ref{thm1} and Theorem \ref{thmdisktype} cannot be dropped. Let 
$$
\Omega = \{(x,y,z) \in \R^3: 1 < \sqrt{x^2+y^2} < 2,\,-1 < z < 1\}.
$$
If $V_k$ is the varifold given by the annulus 
$A = \{(x,y,0):  1 < \sqrt{x^2+y^2} < 2\}$ with multiplicity $k$, 
then $V_k$ has zero curvature and orthogonal boundary in 
$\Omega$, but the mass ${\bf M}(V_k) = 3 k \pi$ goes to infinity.

By a slight modification one also gets a counterexample 
of disk-type surfaces with bounded curvature and orthogonal
boundary in a smooth domain. Namely one can add to the previous example two copies
of the vertical strip $\{(x,0,z): 1 < x < 2,\,0 < z < 1\}$, 
thereby obtaining a surface with orthogonal boundary which 
can be parametrized on a rectangle. A suitable 
approximation yields a smooth example.} 
\end{example}

%\begin{example}{\em In the estimates the second fundamental 
%form cannot be replaced by the mean curvature. For instance 
%let $\Sigma$ be one period of a triply periodic minimal 
%surface in a unit cube. Filling the unit cube with $N^3$ 
%copies translated and scaled by factor $\frac{1}{N}$, one 
%obtains a sequence $\Sigma_N$ of minimal surfaces with 
%orthogonal boundary. One checks that $|\Sigma_N| = N \,|\Sigma|$,
%$|\partial \Sigma_N| \sim N\,|\partial \Sigma|$ and 
%$|\chi(\Sigma_N)| = N^3 |\chi(\Sigma)|$.}
%\end{example}

\section{Minimizers with Orthogonal Boundary}

{\em Proof of Theorem \ref{thmexistence}. }
If $\Omega$ has an orthogonal slice $\Delta$, then it 
is a minimizer with curvature zero and the theorem follows. 
Otherwise let $V_k \in {\bf CV}^2_\perp(\Omega)$ be a minimizing 
sequence. By Theorem \ref{thmvarifolds} we then have
\begin{equation}
\label{equpperbound}
{\bf M}(V_k) + {\bf M}(\Gamma_k) \leq C\,\|B_{V_k}\|_{L^2(V_k)}^2
\quad \mbox{ where } C = C(n,\Omega).
\end{equation}
Lemma \ref{lemmacompactness} yields that after passing to a subsequence, 
$V_k \to V$ and $\Gamma_k \to \Gamma$ as varifolds in $\R^n$, where 
$V$, $\Gamma$ satisfy Definition \ref{defadmissible} and hence
$V \in {\bf CV}^2_\perp(\Omega)$, and further
$$
\|B_V\|_{L^2(V)} \leq \liminf_{k \to \infty} \|B_{V_k}\|_{L^2(V_k)}.
$$
By Remark \ref{rmkfirstvariation}, we have for any 
$\phi \in C^1_c(\R^n,\R^n)$
$$
\delta V_k(\phi) = 
- \int_{\R^n} \langle H_{V_k},\phi \rangle\,d\mu_{V_k}
- \int_S \langle \nu^S(x),\phi(x) \rangle\,d\sigma_{\Gamma_k}. 
$$
Here $H_{V_k}(x) = {\rm tr\,}B_{V_k}(x,T_x \mu_{V_k})$ and 
$\sigma_{\Gamma_k}$ is the weight measure of $\Gamma_k$.
In particular 
$$
|\delta V_k(\phi)| \leq 
\Big(\|H_{V_k}\|_{L^2(V_k)} {\bf M}(V_k)^{\frac{1}{2}} + {\bf M}(\Gamma_k)\Big)
\|\phi\|_{C^0(\R^n)} \leq C\,\|\phi\|_{C^0(\R^n)}.
$$
Now Allard's integral compactness \cite{All72} implies 
$V \in {\bf ICV}^2_\perp(\Omega)$. All that remains 
to show is that $V \neq 0$, i.e. we need a lower bound
\begin{equation}
\label{eqlowermassbound}
\liminf_{k \to \infty} {\bf M}(V_k) > 0.
\end{equation}
This is achieved by the following three lemmas. Namely, 
if we had ${\bf M}(V_k) \to 0$, then Lemma \ref{lemmalowerinterior} 
and Lemma \ref{lemmalowerboundary} would imply
$$
\liminf_{k \to \infty} \|B_{V_k}\|_{L^2(V_k)}^2  \geq 8\pi.
$$
But then $V_k$ is not a minimizing sequence by Lemma 
\ref{lemmacomparison}, a contradiction. \kasten

\begin{lemma} \label{lemmacomparison} For any bounded domain 
$\Omega \subset \R^n$ of class $C^3$ there exists a surface
$\Sigma \in {\bf S}^2_\perp(\Omega)$ with varifold curvature
\begin{equation}
\label{eqcopmparison}
\int_{\Sigma} |B_\Sigma|^2\,d\mu_\Sigma < 8\pi.
\end{equation}
\end{lemma}

{\em Proof. }Let $S \subset \R^n$ be given by a graph representation
$y = u(z)$ where $z \in \R^{n-1}$ with $|z| \leq 1$, such that
$$
u(0) = 0  \quad \mbox{ and } \quad Du(0) = 0.
$$
The upward unit normal of the graph is given by
$$
\nu(z) = \frac{(-Du(z),1)}{\sqrt{1+|Du(z)|^2}}.
$$
We identify $x \in \R^2$ with $(x,0) \in \R^{n-1}$, and consider
the $2$-dimensional halfsphere
$$
S_{+} = \{(x,y): x \in \R^2,\,y \geq 0,\,|x|^2 + y^2 = 1\} \subset \R^n.
$$
We now define the map
$$
f:S_{+} \to \R^n,\,f(x,y) = (x,u(x)) + y\, \nu(x).
$$
For $u(x) \equiv 0$ we have $f(x,y) = (x,0) + y (0,1) = (x,y)$.
Now consider
$$
u_\lambda(z) = \begin{cases}
\frac{1}{\lambda} u(\lambda z) & \mbox{ for }\lambda \neq 0,\\
0 & \mbox{ for } \lambda = 0.
\end{cases}
$$
Clearly $u_{-\lambda}(z) = - u_\lambda(-z)$, and
$Du_\lambda(z) = Du(\lambda z)$. Moreover $u_\lambda$ is
$C^1$ in both variables $(\lambda,z)$, in fact we have
$$
u_\lambda(z) = \int_0^1 \frac{d}{dt} u_\lambda(tz)\,dt
= \int_0^1 Du(\lambda tz)\cdot z\,dt.
$$
The graph of $u_\lambda$ has the upward unit normal
$$
\nu_\lambda(z) = \frac{(-Du(\lambda z),1)}{\sqrt{1+|Du(\lambda z)|^2}}.
$$
We obtain the variation
$$
f_\lambda:S_{+} \to \R^n,\,f_\lambda(x,y) = (x,u_\lambda(x)) + y\,\nu_\lambda(x).
$$
Clearly $u_0(x) = 0$, $\nu_0(x) = (0,1)$ and $f_0(x,y) = (x,y)$. Now
$$
\frac{\partial}{\partial \lambda} Du(\lambda z) \cdot e_i|_{\lambda = 0}
= D^2 u(0)(z,e_i) = h^S(z,e_i).
$$
Here $h^S$ is the second fundamental form of $S$ at $(z,u(z)) = (0,0)$. Next compute
$$
\frac{\partial}{\partial \lambda} u_\lambda(z)|_{\lambda = 0}
= \int_0^1 D^2 u(0)(tz,z)\,dt = \frac{1}{2} h^S(z,z),
$$
$$
\frac{\partial}{\partial \lambda}\nu_\lambda(z)|_{\lambda = 0}
= \frac{\partial}{\partial \lambda}\,
\frac{(-Du_\lambda(z),1)}{(1+|Du_\lambda(z)|^2)^{\frac{1}{2}}}|_{\lambda = 0}
= - \sum_{i=1}^{n-1} h^S(z,e_i)e_i = D\nu(0)z.
$$
In summary the velocity field of the variation is
$$
\phi(x,y): = \frac{\partial}{\partial \lambda}f_\lambda(x,y)|_{\lambda = 0}
= \frac{1}{2}h^S(x,x) (0,1) + y\,D\nu(0)x.
$$
The surface $f(x,y)$ meets the graph of $u(z)$ orthogonally
along $S^1 = \partial S_{+} \subset \R^2$. To see this, let
$x \in S^1$ and choose a unit vector $v \in \R^2$, $v \perp x$.
At the boundary we have one tangent vector of $f(x,y)$ given by
$$
\frac{d}{d\alpha} f(\cos \alpha\, x + \sin \alpha\, v,0)|_{\alpha =0}
= (v,Du(x)v).
$$
A second tangent vector is obtained by
$$
Df(x,0) \cdot (0,1) = \nu(x) = \frac{(-Du(x),1)}{\sqrt{1+|Du(x)|^2}}.
$$
The first vector is tangent to $f|_{S^1}$, the second is
orthogonal to the first and is normalized. Hence the
second vector is the interior co-normal $\eta$ along the
boundary, and is equal to $\nu(x)$, so that $f(x,y)$ 
satisfies the orthogonal boundary constraint. Since the same computation  applies to $f_\lambda$  
the family $f_\lambda(x,y)$ is admissible for comparison.
We have the following first variation formula
for the Willmore functional with boundary terms, see 
Theorem 1 in \cite{AK16}:
$$
\frac{d}{d\lambda} {\cal W}(f_\lambda)|_{\lambda = 0} =
\frac{1}{2} \int_{\Sigma} \langle \vec{W}(f),\phi \rangle\,d\mu_g
+ \frac{1}{2} \int_{\partial \Sigma} \omega(\eta)\,ds_g.
$$
The boundary term is given by
$$
\omega(\eta) = 2 \langle \phi,\nabla_\eta H \rangle
- d \langle \phi,H \rangle(\eta)
- \frac{1}{2} |H|^2 \langle \phi, \eta \rangle.
$$
In the case of $S_{+}$ we have $\vec{W} = 0$, $H(x,y) = - 2 (x,y)$
and $\nabla H = 0$. Thus
$$
\langle \phi,H \rangle =
- h^S(x,x) y - 2y\, \langle D\nu(0) x, x \rangle =
h^S(x,x) y.
$$
Moreover, $\eta = (0,1)$ along $\partial S_{+}$ and thus
$- d \langle \phi,H \rangle (\eta) = - h^S(x,x)$.
For the third term on the boundary, we get directly
$- \frac{1}{2} |H|^2 \langle \phi, \eta \rangle
= - h^S(x,x)$. We conclude that $\omega(\eta) = - 2 h^S(x,x)$ and thus
$$
\frac{d}{d\lambda} {\cal W}(f_\lambda)|_{\lambda = 0} =
- \int_{S^1} h^S(x,x)\,ds(x) = - 2\pi\, {\rm tr}_{\R^2}h^S.
$$
We may now assume that the 2-plane ``$\R^2$'' in the definition of $S_+$ is spanned by the directions of the two largest principal curvatures of $S$. Furthermore we may assume that our graph representation was chosen around a point where the sum of the two largest principal curvatures is positive. (Such point always exists due to the compactness of $\Omega$). In this case we have achieved

%If the sum of the two biggest principal
%curvatures in direction $(0,1)$ is positive, then
$$
\frac{d}{d\lambda} {\cal W}(f_\lambda)|_{\lambda = 0} < 0.
$$
On the other hand, for $\lambda = 0$ the integral of $|A^\circ|^2$
is zero and hence minimal, therefore its derivative vanishes.
Using $|A|^2 = |A^\circ|^2 + \frac{1}{2}|H|^2$ we conclude
$$
\frac{d}{d\lambda} \int_{S_{+}}
|A_{f_\lambda}|^2\,d\mu_{f_\lambda}|_{\lambda = 0} < 0.
$$
Thus for small $\lambda > 0$ we have comparison surfaces $f_\lambda$
with curvature energy below the half sphere. \kasten

\begin{lemma} \label{lemmalowerinterior}
Let $V_k \in {\bf ICV}(\Omega)$ be a sequence with $\mu_{V_k} \to 0$ 
locally in $\Omega$. Assume there are points $x_k \in {\rm spt\,}\mu_{V_k}$ 
such that $x_k \to x_0 \in \Omega$. Then
\begin{equation}
\label{eqlowerinterior}
\liminf_{k \to \infty} \|B_{V_k}\|_{L^2(V_k)}^2 \geq 16\pi 
= \|B_{{\mathbb S}^2}\|_{L^2({\mathbb S}^2)}^2.
\end {equation}
\end{lemma}

{\em Proof.} Let $0 < \varrho < \mathrm{dist}(x_0,\partial \Omega)$. 
By the monotonicity formula, see equations (A.8) and (A.10) 
in \cite{KS04}, we have for sufficiently large $k$
\begin{eqnarray*}
1 & \leq & \frac{\mu_{V_k}(B_\varrho(x_k))}{\pi \varrho^2}
+ \frac{1}{16 \pi} \int_{B_\varrho(x_k)} |H_{V_k}|^2\,d\mu_{V_k}
+ \frac{1}{2\pi \varrho^2} 
\int_{B_\varrho(x_k)} \langle H_{V_k}(x),x-x_0 \rangle\,d\mu_V\\
& \leq &
\Big(\frac{1}{16 \pi} + \ve\Big) \int_{B_\varrho(x_k)} |H_{V_k}|^2\,d\mu_{V_k}
+ C(\ve)\, \frac{\mu_{V_k}(B_\varrho(x_k))}{\varrho^2}.
\end{eqnarray*}
Letting $k \to \infty$ and then $\ve \searrow 0$ we infer
$$
\liminf_{k \to \infty} \int_{B_\varrho(x_k)} |H_{V_k}|^2\,d\mu_{V_k} \geq 16\pi.
$$
From $H_V(x) = {\rm tr\,}B(x,T_x\mu_V)$ we obviously have
$|B_V|^2 \geq \frac{1}{2} |H_V|^2$, but this improves 
by the factor $2$. Namely, Hutchinson showed the following 
relations for any $v \in \R^n$:
\begin{itemize}
\item $B(x,P)(v,P) \subset P^\perp$ and $B(x,P)(v,P^\perp) \subset P$
\quad (\cite{Hut86}, Proposition 5.2.4(iii))
\item $\langle B(v,\tau),\nu \rangle = \langle B(v,\nu), \tau \rangle$
for $\tau \in P$, $\nu \in P^\perp$ \quad
(\cite{Hut86}, Proposition 5.2.4(i)).
\end{itemize}
Moreover, Mantegazza proved in \cite[Theorem 5.4]{Man96} that
$B(x,P)(P^\perp v,\,\cdot\,) \equiv 0$ for $P = T_x \mu$, if 
$V$ is integer rectifiable. Thus 
$H_V(x) = {\rm tr\,}B^\perp(x,T_x \mu)$ and
\begin{equation}
\label{eqneededbymarius}
\int_{\R^n} |H_V|^2 \,d\mu_V
\leq 2 \int_{\R^n} |B_V^\perp|^2\,d\mu_V
= \int_{\R^n} |B_V|^2 \,d\mu_V.
\end{equation}
This finishes the proof of the lemma. \kasten

\begin{lemma}\label{lemmalowerboundary} Let $\Omega \subset \R^n$ be
a bounded domain of class $C^2$, and let $V_k \in {\bf ICV}^2_{\perp}(\Omega)$
be a sequence with ${\bf M}(V_k) \to 0$ and ${\rm spt\,}V_k \subset U_{\delta_k}(S)$
where $\delta_k \to 0$. Then
\begin{equation}
\label{eqboundarybound}
\liminf_{k \to \infty}
\int_{\overline{\Omega}} |B_{V_k}|^2\,dV_k \geq 8\pi
= \int_{{\mathbb S}^2_{+}} |B_{{\mathbb S}^2}|^2\,dV_{{\mathbb S}^2}.
\end{equation}
\end{lemma}

%We apply this for $V = V_k$. For $k$ large, the assumptions
%of statement (c) in Proposition \ref{propreflectedvarifold}
%are satisfied. Using Willmore's inequality \cite{Wil65}, see
%\cite{KS12} for a varifold version, we get
%\begin{eqnarray*}
%16\pi & \leq & \int_{\R^n} |H_{W_k}|^2\,d\mu_{W_k}\\
%& \leq & 2 \int_{\overline{\Omega}} |H_{V_k}|^2\,d\mu_{V_k}
%+ C \ve \int_{\overline{\Omega}} |B_{V_k}|^2\,d\mu_{V_k}
%+ C(\Lambda,\ve) {\bf M}(V_k)\\
%& \leq & (2 + C\ve)
%\int_{\overline{\Omega}} |B_{V_k}|^2\,d\mu_{V_k}
%+ C(\Lambda,\ve) {\bf M}(V_k).
%\end{eqnarray*}
%The lemma follows by letting $k \to \infty$ and then
%$\ve \searrow 0$.
%\kasten

The rest of this section is devoted to the proof of Lemma 
\ref{lemmalowerboundary}. A monotonicity formula for 
$m$-dimensional varifolds with free boundary and mean 
curvature in $L^p$, $p > m$, was proved by Gr\"uter
and Jost in \cite{GJ86}, see also De Masi \cite{DeM20}.  
However, it seems unclear whether the argument of 
\cite{GJ86} extends to the critical case $p = m = 2$.
Our strategy is to reflect the varifold across $S$. 
This approach is not possible in the setting of \cite{GJ86}, 
because the estimate for the mean curvature of the 
reflected varifold involves the second fundamental form, 
see equation (\ref{eqreflectedvarifoldenergy}) below.\\
\\
In the following, it will be convenient to 
pullback the geometry by the reflection map.

\begin{definition} Let $G = (g_{ij}) \in C^1(U,\R^{n \times n})$ be a Riemannian
metric on the open set $U \subset \R^n$. For an $m$-varifold $V$ on $U$
we consider the functional
\begin{equation}
\label{eqriemannianvarifold}
V^g(\phi) = \int_{U \times G(m,n)} \phi(x,P)\, JG(x,P)\,dV(x,P) \quad
\mbox{ for }\phi \in  C^0_c(U \times G(m,n)).
\end{equation}
Here $JG(x,P)$ denotes the Riemannian volume element. With respect to a
standard orthonormal basis $\tau_\alpha$ of $P$ where $\alpha = 1,\ldots,m$,
we have
\begin{equation}
JG(x,P) = \big(\det g(x)(\tau_\alpha,\tau_\beta) \big)^{\frac{1}{2}}
= \big(\det PG(x)|_P\big)^{\frac{1}{2}}.
\end{equation}
The Riemannian mass of $V$ is given by
\begin{equation}
\label{eqriemannianmass}
M^g(V) = \int_{U \times G(m,n)} JG(x,P)\,dV(x,P).
\end{equation}
The first variation of $V^g$ with respect to a vector field
$X \in C^1_c(U,\R^n)$ is defined by
\begin{equation}
\label{riemannfirstvariation}
\delta V^g(X) = \frac{d}{dt}  M^g((\varphi_t)_\sharp V)|_{t=0},
\end{equation}
where $\varphi_t$ is the flow associated to $X$. 
\end{definition}

\begin{lemma}[Riemannian first variation] \label{lemmariemannfirstvariation}
Let $X \in C^1_c(U,\R^n)$ be a vector field with associated flow
$\varphi_t$. Assume that $V$ is an integer rectifiable $m$-varifold in $U$
with weak second fundamental form $B \in L^1(V)$. Then we have the first
variation formula
\begin{equation}
\label{eqriemannfirstvariation}
\delta V^g(X)
= - \int_U \langle G(x) H_g(x),X(x) \rangle\,JG(x,P(x))\,d\mu_V(x),
%\label{eqriemannfirstvariation}
\end{equation}
where $\mu_V$ is the weight measure of $V$ and $P(x) = T_xV$ is the
approximate tangent space. The vector $H_g$ decomposes as
$H_g(x) = \sum_{i=1}^3 \phi_i(x,P(x))$ where
\begin{eqnarray}
\label{eqphi1}
\phi_1(x,P)^l & = & \frac{1}{JG(x,P)}\, g^{lm}(x) \frac{\partial}{\partial P^j_k}
\big(JG\,Q^i_m\big)(x,P)\, B^k_{ij}(x,P),\\
\label{eqphi2}
\phi_2(x,P) & = & \frac{1}{JG(x,P)}\, G(x)^{-1}\frac{\partial}{\partial x^j}
\big(JG\,\,Q\big)(x,P)^{\rm T}Pe_j,\\
\label{eqphi3}
\phi_3(x,P)^i & = &
- \frac{1}{2} g(x)^{ij} \big\langle Q(x,P) G(x)^{-1} \partial_{x^j} G(x),P \big\rangle.
\end{eqnarray}
Here $Q(x,P)$ is the projection onto $P$ with respect to the inner product $g(x)$.
\end{lemma}

{\em Proof. }For a diffeomorphism $\varphi \in C^1(U,U)$ the pushforward varifold
is defined by
\begin{equation}
\label{eqpushforward}
\varphi_\sharp V (\psi) =
\int_{U \times G(m,n)} \psi(\varphi(x),D\varphi(x)P)\,J\varphi(x,P)\,dV(x,P).
\end{equation}
Here by $D\varphi(x)P$ we mean the image subspace, and $J\varphi(x,P)$ is the
standard Jacobian of $\varphi$ on $P$. The Riemannian mass of $\varphi_\sharp V$
is then
\begin{eqnarray*}
M^g(\varphi_\sharp V) & = & \int_{U \times G(m,n)} JG(y,Q)\,d(\varphi_\sharp V)(y,Q)\\
& = & \int_{U \times G(m,n)}
JG(\varphi(x),D\varphi(x)P)\,J\varphi(x,P)\,dV(x,P).
\end{eqnarray*}
Now consider the identity maps $I:(P,g(x)) \to P$ and
$J:D\varphi(x)P \to \big(D\varphi(x)P,g(\varphi(x))\big)$. Let
$A(x): P \to D\varphi(x)P$, $A(x)v = D\varphi(x)v$. We have
\begin{eqnarray*}
J^g \varphi(x,P) &  = &
\det \big((J A(x) I)^\ast\, J A(x) I\big)^{\frac{1}{2}}\\
& = & \det (J^\ast J)^{\frac{1}{2}} \,\det(A(x)^\ast A(x))^{\frac{1}{2}}\,
\det(I^\ast I)^{\frac{1}{2}}\\
& = &  \frac{JG(\varphi(x),D\varphi(x)P)}{JG(x,P)}{J\varphi(x,P)}.
\end{eqnarray*}
It follows that
$$
M^g(\varphi_\sharp V) = \int_{U \times G(m,n)} J^g \varphi(x,P)\,dV^g(x,P)
\quad \mbox{ where }dV^g = JG\,\,dV.
$$
Now assume that $\varphi = \varphi_t$ is the flow of a vector field
$X \in C^1_c(U,\R^n)$. Choose an orthonormal basis $v_1,\ldots,v_m$
of $P$ with respect to $g(x)$, and calculate
\begin{eqnarray*}
\frac{\partial}{\partial t} J^g\varphi_t(x,P)|_{t=0} & = &
\frac{\partial}{\partial t} \Big(\det g(\varphi_t(x))
(D\varphi_t(x) v_\alpha, D\varphi_t(x) v_\beta)\Big)^{\frac{1}{2}}|_{t=0}\\
& = & g(x)(DX(x)v_\alpha,v_\alpha) + \frac{1}{2} D_X g(x)(v_\alpha,v_\alpha)\\
& = & {\rm tr\,}\big(Q(x,P) DX(x)|_P\big)
+ \frac{1}{2} {\rm tr\,}\big(Q(x,P) G(x)^{-1} D_X G(x)|_P\big)\\
& = & \langle Q(x,P)\,DX(x),P \rangle
+ \frac{1}{2} \langle Q(x,P)\,G(x)^{-1}D_X G(x),P \rangle.
\end{eqnarray*}
We write
\begin{eqnarray*}
JG(x,p) \langle Q(x,P) DX(x), P \rangle & = &
\langle D_x \big(JG(x,P)Q(x,P) X(x)\big), P \rangle\\
&& - \sum_{j=1}^n \big\langle
\frac{\partial}{\partial x^j} \big(JG(x,P)Q(x,P)\big)\,X(x),Pe_j \big\rangle.
\end{eqnarray*}
We now apply the curvature varifold property with the test
vector field $Y(x,P) = JG(x,P) Q(x,P) X(x)$. Noting that
$Y$ is $C^1$ and has compact support, we get
\begin{eqnarray}
\label{eqpartialintegration}
&& \int_{U \times G(m,n)} \langle D_x \big(JG(x,P)Q(x,P) X(x)\big), P \rangle\,dV(x,P)\\
\nonumber
& = & - \int_{U \times G(m,n)} \big(D_P Y \cdot B + \langle {\rm tr\,}B,Y \rangle\big)\,dV(x,P).
%&& - \int_{G_{m-1}(TS)} \big\langle \nu^S(x),Y(x,\nu^S(x) \wedge Q)\big\rangle \,d\Gamma(x,Q).
\end{eqnarray}
Since $V$ is integer rectifiable, it has weak mean curvature
$$
H^k_V(x) = \sum_{i=1}^n B_{ii}^k(x,P(x)) \in L^1(\mu_V),
$$
where $P(x)$ is the projection onto $T_x\mu$
and $\mu_V$ is the weight measure. Now we use 
$H(x) \perp T_x\mu$ $\mu$-almost-everywhere,
see Brakke \cite[Chapter 5]{Bra78} or
Mantegazza \cite[Theorem 5.4]{Man96}. As $Y(x,P(x)) \in T_x\mu$, 
the last term in (\ref{eqpartialintegration}) vanishes. We conclude
\begin{eqnarray*}
\delta V^g(X) & = & \frac{d}{dt} M^g((\varphi_t)_\sharp V)|_{t=0}\\
& = & - \int_{U \times G(m,n)} D_P Y \cdot B \,dV(x,P)\\
&& - \int_U \sum_{j=1}^n \big\langle \partial_{x^j} (JG\, \,Q)(x,P)X(x),Pe_j \big\rangle\,dV(x,P)\\
&& + \frac{1}{2} \int_{U \times G(m,n)}  \langle Q(x,P) G(x)^{-1} D_X G(x),P \rangle\,
JG(x,P)\,dV(x,P).
\end{eqnarray*}
%&&  - \int_{G_{m-1}(TS)} \big\langle \nu^S(x),Y(x,\nu^S(x) \wedge Q)\big\rangle \,d\Gamma(x,Q).
%\end{eqnarray*}
Finally, we rewrite this in Riemannian form, that is
$$
\delta V^g(X)
= - \int_U \langle G(x) H_g(x),X(x) \rangle\,JG(x,P(x))\,d\mu_V(x).
$$
Now one readily checks this formula for $H_g$ as given in
the statement of the Lemma. \kasten

\begin{lemma} [local energy bound] \label{lemmalocalenergy} Assume that
in addition to the assumptions of the previous lemma, we have
\begin{equation}
\label{eqlocalassumptions}
|G(x) - \mathrm{E}_n| \leq \ve \quad \mbox{ and } \quad |DG(x)| \leq \Lambda.
\end{equation}
Then we can estimate, for $P = P(x)$,
\begin{equation}
\label{eqlocalenergy}
\|H_g(x)\|^2_g \,JG(x,P) \leq |H(x)|^2 + C\ve\,|B(x,P)|^2 + C(\ve,\Lambda).
\end{equation}
\end{lemma}

{\em Proof. } We consider the sets
\begin{eqnarray*}
{\cal G} &  = &
\{G \in \R^{n \times n}: G^{\rm T} = G,\,|G-{\rm E}_n| \leq \frac{1}{2}\},\\
{\cal M} & = & \{M \in L(\R^n,\R^{n \times n}): |M| \leq \Lambda\}.
\end{eqnarray*}
For given $G \in {\cal G}$ and $P \in G(m,n)$, we choose a standard orthonormal
basis $\tau_1,\ldots,\tau_m$ of $P$ and define $g_{\alpha \beta} =
\langle G \tau_\alpha,\tau_\beta \rangle$. The following quantities are
well-defined:
\begin{eqnarray*}
J(G,P) & = & \big(\det g_{\alpha \beta} \big)^{\frac{1}{2}},\\
Q(G,P)e_j & = &  g^{\alpha \beta} \langle G \tau_\beta,e_j \rangle\,\tau_\alpha.
\end{eqnarray*}
To see that the functions are smooth, one can choose a smooth orthonormal
basis $\tau_\alpha(P)$ for $P$ near a point $P_0 \in G(m,n)$. Such a basis
can be obtained by taking a basis $\tau_\alpha$ for $P_0$ and applying
the Gram-Schmidt process to the vectors $P \tau_\alpha$.
Putting $G = G(x)$ and $M = DG(x)$, one verifies that the three
functions $\phi_i(x,P)$ for $P = P(x)$ have the following structure:
\begin{eqnarray*}
\phi_1(x,P) & = & f_1(G,P)^{ij}_k \,B_{ij}^k,\\
\phi_2(x,P) & = & f_2(G,P) \ast M,\\
\phi_3(x,P) & = & f_3(G,P) \ast M.
\end{eqnarray*}
(Here, $f(G,P) \ast M$ refers to a sum of products of terms of expressions which are smooth in $G,P$ with entries of $M$.) 
For $G = {\rm E}_n$ and $M = 0$ we compute precisely
$$
\phi_1^l = \delta^{lm} \frac{\partial P^i_m}{\partial P^j_k} B_{ij}^k
= B^l_{ii} = H^l \quad \mbox{ and } \quad
\phi_2 = \phi_3 = 0.
$$
It follows that
\begin{eqnarray*}
|H_g(x) - H(x)| & \leq &  C\,|f_1(G,P)-f_1({\rm E}_n,P)|\,|B|
+ C\,\big(|f_2(G,P)|+|f_3(G,P)|\big)\, |M|\\
& \leq & C \ve \,|B| + C \Lambda.
\end{eqnarray*}
We can now estimate further
\begin{eqnarray*}
\big||H_g|^2-|H|^2\big| & \leq &
\big(2|H|+ C \ve \,|B| + C \Lambda\big) \cdot \big(C \ve |B| + C \Lambda)\\
& \leq & C (|B| + C \Lambda)(C \ve |B| + C \Lambda)\\
& \leq & C \ve |B|^2 + C\Lambda |B| + C \Lambda^2\\
& \leq & C \ve |B|^2 + C(\ve) \Lambda^2.
\end{eqnarray*}
Finally we use $\langle Gv,v \rangle \leq (1+C\ve)|v|^2$ and conclude
\begin{eqnarray*}
\|H_g\|_g^2\,JG & \leq & (1+C\ve) |H_g|^2\\
& \leq & (1+C\ve)(|H|^2 + C \ve |B|^2 + C(\ve) \Lambda^2)\\
& \leq & |H|^2 + C \ve |B|^2 + C(\ve) \Lambda^2.
\end{eqnarray*}
The lemma is proved. \kasten

We now turn to the problem of mapping a varifold by a diffeomorphism.

\begin{lemma} Let $\sigma \in C^2(U,\hat{U})$ be a diffeomorphism
of open sets $U,\hat{U} \subset \R^n$, with pullback metric
$g_{ij} = \langle \partial_i\sigma,\partial_j \sigma \rangle$.
Let $V$ be an integer rectifiable varifold with compact
support in $U$, having second fundamental form $B \in L^2(V)$.
Then the pushforward varifold $\hat{V} = \sigma_\sharp V$ has
weak mean curvature $\hat{H} \in L^2(\hat{\mu})$, 
for $\hat{\mu} = \mu_{\hat{V}}$, given by
\begin{equation}
\label{eqdiffeomean}
\hat{H} \circ \sigma = D\sigma \cdot H_g.
\end{equation}
Moreover, we have for $JG(x) = JG(x,T_x \mu_V)$ and any Borel set $E \subset U$
\begin{equation}
\label{eqdiffeonorm}
\int_{\sigma(E)} |\hat{H}(y)|^2\,d\hat{\mu}(y) =
\int_E \|H_g(x)\|_g^2\,JG(x)\,d\mu(x).
\end{equation}
\end{lemma}

{\em Proof. }By definition of the varifold pushforward, we have
$$
\sigma_\sharp V(\eta) =
\int_U \eta(\sigma(x),D\sigma(x)P) J\sigma(x,P)\,dV(x,P).
$$
For an orthonormal basis $\tau_\alpha$ of $P$ we have
$\langle D\sigma(x) \tau_\alpha,D\sigma(x) \tau_\beta \rangle
= g(x)(\tau_\alpha,\tau_\beta)$, thus
$$
{\bf M}(\sigma_\sharp V) = \int_U J\sigma(x,P)\,dV(x,P) =
\int_U JG(x,P)\,dV(x,P) = M^g(V).
$$
Now for given $Y \in C^1_c(\hat{U},\R^n)$ we let
$X \in C^1_c(U,\R^n)$ be the $\sigma$-related vector field,
i.e. $Y \circ \sigma = D\sigma \cdot X$. Then we have
$\psi_t \circ \sigma = \sigma \circ \varphi_t$ for the
related flows, and
$$
{\bf M}((\psi_t)_\sharp\, \sigma_\sharp V) = {\bf M}((\psi_t \circ \sigma)_\sharp V)
= {\bf M}((\sigma \circ \varphi_t)_\sharp V)
= {\bf M}(\sigma_\sharp\, (\varphi_t)_\sharp V)
= M^g ((\varphi_t)_\sharp V).
$$
Therefore we obtain
\begin{eqnarray*}
\delta \hat{V}(Y) & = &
\frac{d}{dt} {\bf M}((\psi_t)_\sharp (\sigma_\sharp V))|_{t=0}\\
& = & \frac{d}{dt} M^g((\varphi_t)_\sharp V)|_{t=0}\\
& = & - \int_U \langle G(x) H_g(x),X(x)\rangle\,JG(x,P(x))\,d\mu(x).
%& & - \int_{G_{m-1}(TS)} \langle \nu^S(x),X(x)\rangle\,d\Gamma(x,Q).
\end{eqnarray*}
Transforming back to $\sigma_\sharp V$ yields
\begin{eqnarray*}
&&\int_U \langle G(x) H_g(x),X(x)\rangle\,JG(x,P(x))\,d\mu(x)\\
& = & \int_{U \times G(m,n)} \langle D\sigma(x) H_g(x), D\sigma(x) X(x)\rangle
J\sigma(x,P)\,dV(x,P)\\
& = & \int_{U \times G(m,n)} \langle D\sigma(\sigma^{-1}(y)) H_g(\sigma^{-1}(y)),
Y(y) \rangle|_{y = \sigma(x)}\, J\sigma(x,P)\,dV(x,P)\\
& = & \int_{\hat{U} \times \R^n}
\big\langle D\sigma(\sigma^{-1}(y)) H_g(\sigma^{-1}(y)), Y(y) \big\rangle\,
d(\sigma_\sharp V)(y,P).
\end{eqnarray*}
%For the boundary term, we use the properties of $\sigma$ to get
%\begin{eqnarray*}
%\langle \nu^S(x),X(x)\rangle
%& = & \langle \nu^S(x),D\sigma(x)^{-1}Y(\sigma(x)) \rangle\\
%& = & \langle D\sigma(x) \nu^S(x), Y(\sigma(x)) \rangle\\
%& = & - \langle \nu^S(x), Y(x) \rangle.
%\end{eqnarray*}
We conclude
\begin{equation}
\label{eqvariationreflect}
\delta \hat{V}(Y) =
- \int_{\hat{U}} \langle \hat{H}, Y \rangle\,d\hat{\mu}
\quad \mbox{ where } \hat{H}(y) =  D\sigma(x) H_g(x)|_{x = \sigma^{-1}(y)}.
%+ \int_{G_{m-1}(TS)} \langle \nu^S(y),Y(y)\rangle\,d\Gamma(y,Q),
\end{equation}
Substituting $X(x) = H_g(x)$ in the previous calculation, we see finally
$$
\int_U \|H_g(x)\|_g^2 JG(x)\,d\mu(x) =
\int_{\hat{U}} |\hat{H}(y)|^2\,d\hat{\mu}(y).
$$
\kasten

 We now specialize to the situation where $\Omega \subset \R^n$ is a
bounded $C^2$ domain with boundary $S = \partial \Omega$, and
consider  the reflection
\begin{equation}
\label{eqreflection}
\sigma:U \to U,\,\sigma(x+r\nu^S(x)) = x-r\nu^S(x) \quad
\mbox{ where }(x,r) \in S \times (-\delta,\delta).
\end{equation}
Here $\delta> 0$ is chosen so small that the two-sided tubular neighborhood $U = U_\delta(S) = \{ x + r \nu^S(x) : (x,r) \in S \times (-\delta,\delta) \}$ is  open. Notice further that $\sigma(U) = U$. We put
$g_{ij}(x) = \langle \partial_i \sigma(x),  \partial_j \sigma(x) \rangle$ for all $x \in U$.

\begin{lemma} \label{lemmametriconS} For $x \in S$, $v \in \R^n$ we have,
where $P_S$ is the projection onto $TS$, 
\begin{equation}
\label{eametriconS}
g_{ij}(x) = \delta_{ij} \quad \mbox{ and } \quad
Dg_{ij}(x)v = 4 \langle \nu^S(x),v \rangle \,h^S(x)(P_S e_i,P_S e_j).
\end{equation}
\end{lemma}

{\em Proof. }First the $r$-derivative of (\ref{eqreflection}) gives
$$
D\sigma(x+r \nu^S(x)) \nu^S(x) = - \nu^S(x).
$$
Second let $\gamma(s)$ be a curve in $S$ with $\gamma(0) = x$ and
$\gamma'(0) = v \in T_xS$. Then we compute, denoting by $W(x)$ the
Weingarten map of $S$,
\begin{eqnarray*}
({\rm Id} - rW(x)) v & = & \frac{d}{ds} \big(\gamma(s) - r \nu(\gamma(s))\big)|_{s=0}\\
& = & \frac{d}{ds} \sigma\big(\gamma(s) + r \nu^S(\gamma(s))\big)|_{s=0}\\
& = & D\sigma(x+r\nu^S(x)) \big({\rm Id} + r W(x)\big) v.
\end{eqnarray*}
This implies
$$
D\sigma(x+r\nu^S(x)) v = ({\rm Id}-r W(x))\,({\rm Id}+r W(x))^{-1} v
\quad \mbox{ for }v \in T_xS.
$$
For $r = 0$ we get as expected
$$
D\sigma(x) \nu^S(x) = - \nu^S(x) \quad \mbox{ and } \quad
D\sigma(x) v = v.
$$
This implies $G(x) = D\sigma(x)^{\rm T} D\sigma(x) = {\rm E_n}$
for $x \in S$. Moreover it follows that
$$
DG(x)v = 0 \quad \mbox{ for }x \in S,\,v \in T_x S.
$$
We compute further
$$
\frac{\partial}{\partial r} D\sigma(x+r\nu^S(x)) \nu^S(x)|_{r=0}
= \frac{\partial}{\partial r} (-\nu^S(x)) = 0.
$$
For the tangential derivative we get, again for $v \in T_x S$,
\begin{eqnarray*}
\frac{\partial}{\partial r} D\sigma(x+r\nu^S(x))v|_{r=0}
& = &
\frac{\partial}{\partial r} ({\rm Id}-r W(x))\,({\rm Id}+r W(x))^{-1} v|_{r=0}\\
& = & - 2 W(x) v.
\end{eqnarray*}
Thus we have, again for $v,w \in T_xS$.
$$
\begin{array}{lcl}
\frac{\partial}{\partial r} g(x+r\nu^S(x))(\nu^S(x),\nu^S(x))|_{r=0} & = & 0,\\
\frac{\partial}{\partial r} g(x+r\nu^S(x))(\nu^S(x),v)|_{r=0} & = & 0,\\
\frac{\partial}{\partial r} g(x+r\nu^S(x))(v,w)|_{r=0} & = &
4 h^S(x)(v,w).
\end{array}
$$
The lemma is proved. \kasten

We now compute the vector $H_g(x)$ as
defined in Lemma \ref{lemmariemannfirstvariation} for points $x \in S$.

\begin{lemma} \label{lemmaboundarymeancurvature} For $x \in S$ and
$P = P(x)$ we have the formula
\begin{equation}
\label{eqboundarymeancurvature}
H_g(x) = H(x) - 2\, \Delta_P d_S(x)\, P^\perp \nu^S(x)
+ 4\,P^\perp D\nabla d_S(x) \cdot P\nu^S.
\end{equation}
Here the notation $\Delta_P d_S(x) = {\rm tr}_P D^2 d_S(x)$ is used.
\end{lemma}

{\em Proof. }Let $P \in G(m,n)$ be any fixed subspace, with
Euclidean orthonormal basis $\tau_\alpha$ for $1 \leq \alpha \leq m$.
Put $g_{\alpha \beta}(x) = g(x)(\tau_\alpha,\tau_\beta)$. At a point
$x \in S$ we have
\begin{eqnarray*}
g_{\alpha \beta}(x) & = &
\langle \tau_\alpha,\tau_\beta \rangle = \delta_{\alpha \beta},\\
\partial_{x^k} g_{\alpha \beta}(x) & = &
4 \langle \nu^S(x),e_k \rangle\,D^2 d_S(x)(\tau_\alpha,\tau_\beta).
\end{eqnarray*}
%For convenience, let us introduce the notation
We then have
\begin{equation}
\label{eqderivativejacobi}
\partial_{x^k} JG(x,P) =
2 \langle \nu^S(x),e_k \rangle\,D^2 d_S(\tau_\alpha,\tau_\alpha)
= 2 \langle \nu^S(x),e_k \rangle\,\Delta_P d_S(x).
\end{equation}
Next we compute for $G(x) = (g_{\alpha \beta}(x))$,
\begin{equation}
\label{eqderivativeGP}
\langle \partial_{x^k} G, P \rangle =
4\, \langle \nu^S(x),e_k \rangle\, h^S(x)(\tau_\alpha,\tau_\alpha)
= 4\, \langle \nu^S(x),e_k \rangle \,\Delta_P d_S(x).
\end{equation}
Now we address the $x^k$ derivative of $Q(x,P)$. Using
again $g_{\alpha \beta} = g(\tau_\alpha, \tau_\beta)$ we have
$$
Q(x,P) v = g^{\alpha \beta} g(\tau_\beta,v) \tau_\alpha \quad
\mbox{ for any }v \in \R^n.
$$
Now we have, again for any $v \in \R^n$,
\begin{eqnarray*}
\partial_{x^k} g^{\alpha \beta} & = &
- 4 \langle \nu^S(x),e_k \rangle\,D^2 d_S(x)(\tau_\alpha,\tau_\beta),\\
\partial_{x^k} g(\tau,v) & = & 4 \langle \nu^S(x),e_k \rangle\,
D^2 d_S(x)(\tau,v).
\end{eqnarray*}
Inserting we have
$$
\partial_{x^k} Q(x,P) v =
- 4 \langle \nu^S(x),e_k \rangle\,D^2 d_S(x)(\tau_\alpha,\tau_\beta)
\langle \tau_\beta,v \rangle \,\tau_\alpha
+\, 4\, \langle \nu^S(x),e_k \rangle\,
D^2 d_S(\tau_\alpha,v) \tau_\alpha,
$$
which yields
\begin{eqnarray}
\label{eqderivativeprojection}
\partial_{x^k} Q(x,P) v & = & 4\, \langle \nu^S(x),e_k \rangle\,
D^2 d_S(x)(\tau_\alpha,P^\perp v) \tau_\alpha\\
\nonumber
& = &  4\, \langle \nu^S(x),e_k \rangle\,P\, D\nabla d_S(x) P^\perp v.
\end{eqnarray}
Now $H_g(x) = \sum_{i=1}^3 \phi_i(x,P(x))$,
where for $g_{ij}(x) = \delta_{ij}$ the $\phi_i$ are as follows: 
\begin{eqnarray*}
\phi_1(x,P)^l & = & \partial_{P^j_k}(JG\,Q^i_l)(x,P) B_{ij}^k(x,P),\\
\phi_2(x,P) & = & (\partial_{x^k} JG(x,P))\,Pe_k + \partial_{x^k} Q(x,P)^{\rm T} Pe_k,\\
\phi_3(x,P)^i & = & - \tfrac{1}{2} \langle \partial_i G(x),P \rangle.
\end{eqnarray*}
For $\phi_3$ we have directly
$$
\phi_3(x,P) = - 2 \langle \nu^S(x),e_k \rangle \Delta_P d_S(x) e_k
= - 2 \Delta_P d_S(x) \,\nu^S(x).
$$
For $\phi_2(x,P)$ we first calculate
\begin{eqnarray*}
\langle \partial_{x^k} Q(x,P)^{\rm T} v, w \rangle & = &
\langle v, \partial_{x^k} Q(x,P) w \rangle\\
& = & 4 \langle \nu^S(x),e_k \rangle\,
\langle v, D^2 d_S(\tau_\alpha,P^\perp w) \tau_\alpha \rangle\\
& = &  4 \langle \nu^S(x),e_k \rangle \,
D^2 d_S(x)(Pv,P^\perp w).
\end{eqnarray*}
Inserting we compute further
\begin{eqnarray*}
\phi_2(x,P) & = & 2\, \langle \nu^S(x),e_k \rangle\, \Delta_P d_S(x)\, Pe_k\\
&& + 4\,\langle \nu^S(x),e_k \rangle\,D^2 d_S(x)(Pe_k,P^\perp e_l)\, e_l\\
& = & 2\, \Delta_P d_S(x) P \nu^S(x) + 4\, D^2 d_S(x)(P \nu^S(x),P^\perp e_l) e_l.
\end{eqnarray*}
To simplify further we use
$$
D^2 d_S(x)(P \nu^S(x),P^\perp e_l)e_l
= \langle D\nabla d_S(x) P \nu^S(x),P^\perp e_l \rangle e_l
= P^\perp D\nabla d_S(x) P \nu^S(x).
$$
Thus we arrive at
$$
\phi_2(x,P) =   2\, \Delta_P d_S(x) P \nu^S(x) + 4\,P^\perp D\nabla d_S(x) P \nu^S(x).
$$
It remains to compute $\phi_1(x,P)$. This involves only a partial derivative
with respect to the $P$ variable, hence we can assume that $x \in S$ is
a fixed point, in particular we have for all $P$ that $JG(x,P) = 1$ and
$Q(x,P) = P$. We then get easily
$$
\phi_1(x,P)^l = \frac{\partial P^i_l}{\partial P^j_k} B^k_{ij}(x,P)
= \delta_{ij} \delta_{kl} B^k_{ij}(x,P)
= B^l_{ii}.
$$
In other words we just have
$$
\phi_1(x,P(x)) = H(x).
$$
Collecting the three terms we finally have
\begin{equation}
\label{eqmeancurvatureg}
H_g(x) = H(x) - 2 \Delta_P d_S(x) P^\perp \nu^S(x)
+ 4\,P^\perp D\nabla d_S(x) P \nu^S(x).
\end{equation}
\kasten

\begin{proposition}[reflecting the varifold] \label{propreflectedvarifold}
Let $\Omega \subset \R^n$ be a bounded $C^3$ domain, and let
$V \in {\bf ICV}^m_\perp(\Omega)$ be given with ${\rm spt\,}\mu_V \subset U$.
Here we suppose that $U=U_\delta(S)$ is a tubular neighborhood
of $S = \partial \Omega$ such that the reflection $\sigma \in C^2(U,U)$
across $S$ is defined. Define
\begin{equation}
\label{eqreflectedvarifold}
W(\phi) = V(\phi) + \sigma_\sharp V(\phi) \quad \mbox{ for }
\phi \in C^0_c(\R^n \times G(m,n)).
\end{equation}
Then the integer rectifiable varifold $W$ has the following properties:\\
\\
{\rm (a)} $\mu_W = 2\mu_V$ on $S$, and $T_x\mu_W = T_x \mu_V \subset T_xS$ 
whenever $T_x \mu_V$, $x \in S$, exists.\\
\\
{\rm (b)} $W$ has weak mean curvature $H_W \in L^2(\mu_W)$,
given by
\begin{equation}
\label{eqWmeancurvature}
H_W(x) = \begin{cases}
H_V(x) & \mbox{ for }x \in \Omega,\\
(D\sigma \cdot H_g)(\sigma(x)) &
\mbox{ for }x \in \R^n \backslash \overline{\Omega},\\
\frac{1}{2}\big(H_V(x) + D\sigma(x) H_g(x)\big)
& \mbox{ for } x \in S.
\end{cases}
\end{equation}
Here $H_g$ is as in Lemma \ref{lemmariemannfirstvariation}.
For $x \in S$ we have precisely
\begin{equation}
\label{eqWmeanboundary}
H_W(x) = P^{\top}_S H_V(x) + \Delta_P d_S(x)\, \nu^S(x).
\end{equation}
{\rm (c)} Assume $|G-{\rm E_n}| \leq \ve$ and $|DG| \leq \Lambda$
on $U$. Then $W$ satisfies
\begin{equation}
\label{eqreflectedvarifoldenergy}
\int_{\R^n} |H_W|^2 d\mu_W \leq
2\, \int_{\overline{\Omega}} |H_{V}|^2 d\mu_{V}
+ C \ve \int_{\overline{\Omega}} |B_{V}|^2 d\mu_{V}
+ C(\Lambda,\ve) {\bf M}(V).
\end{equation}
%\item[{\rm (d)}] $W$ has the lower energy bound
%\begin{equation}
%\label{eqreflectedvarifoldwillmore}
%\int_{\R^n} |H_W|^2\,d\mu_W \geq 16\pi.
%\end{equation}
\end{proposition}

{\em Proof. } Assuming $V = {\bf v}(M,\theta)$ we have
$\sigma_\sharp V = {\bf v}(\sigma(M),\theta \circ \sigma)$.
Thus
$$
W = {\bf v}(M \cup \sigma(M), \theta + \theta \circ \sigma)
$$
is integer rectifiable and $\mu_W = 2 \mu_V$ on $S$.
Now if $T_x \mu_V$ exists for $x \in S$, then $T_x \mu_V$
is contained in the halfspace $\langle \nu^S(x),v \rangle \geq 0$,
and hence $T_x\mu_V \subset T_x S$. In general
$T_{\sigma(x)} \mu_{\sigma_\sharp V} = D\sigma(x)T_x \mu_V$,
and hence $T_x \mu_{\sigma_\sharp V} = T_x V = T_x W$
for $x \in S$.\\
\\
We now turn to {\rm (b)}. In the definition of varifold with orthogonal
boundary, we use a test function $\phi \in C^1_c(\R^n)$ to get the first
variation formula
$$
\delta V(\phi)
= - \int_{\overline{\Omega}} \langle H(x),\phi(x) \rangle\,d\mu_V(x)\\
- \int_S \langle \nu^S(x),\phi(x)\rangle\,d\mu_\Gamma(x).
$$
Here $H(x) = \sum_{i=1}^n B_{ii}(x,T_x \mu_V)$ and we used
that $\Gamma^x(G_{m-1}(T_xS)) = 1$ by definition. We now repeat
the computation of Lemma \ref{lemmariemannfirstvariation}.
Since here $\phi$ is not compactly supported in $\Omega$, we
obtain additional boundary terms:
\begin{eqnarray*}
\delta V^g(\phi) & = & - \int_{\overline{\Omega}}
\langle G(x) H_g(x),\phi(x)\rangle \,JG(x,P(x))\,d\mu_V(x)\\
&& - \int_S \langle \nu^S(x), \phi(x) \rangle\,\,d\mu_{\Gamma}(x)
- \int_S \langle PH(x),\phi(x) \rangle\,d\mu_V(x).
\end{eqnarray*}
The last term will eventually be cancelled, however we do not know
at present that the relation $H_V \perp T\mu_V$ holds on $S$.
Replacing $\phi$ by the $\sigma$-related field $\psi$, i.e.
$\phi(\sigma(x)) = D\sigma(x) \psi(x)$, we get
\begin{eqnarray*}
\delta V^g(\psi) & = & - \int_{\overline{\Omega}}
\langle G(x) H_g(x),\psi(x)\rangle \,JG(x,P(x))\,d\mu_V(x)\\
&& + \int_S \langle \nu^S(x), \phi(x) \rangle\,\,d\mu_{\Gamma}(x)
- \int_S \langle PH(x),\phi(x) \rangle\,d\mu_V(x)\\
& = & -  \int_{\R^n \backslash \Omega} \langle
D\sigma(\sigma(y)) \cdot H_g(\sigma(y)),\phi(y)\rangle\,d\mu_{\sigma_\sharp V}(y) \\
&& + \int_S \langle \nu^S(x), \phi(x) \rangle\,\,d\mu_{\Gamma}(x)
- \int_S \langle PH(x),\phi(x) \rangle\,d\mu_V(x).
\end{eqnarray*}
Now $\delta (\sigma_\sharp V)(\phi) = \delta V^g(\psi)$. When adding
the variations of $V$ and $\sigma_\sharp V$, the boundary terms
involving $\Gamma$ cancel, and we obtain
$$
\delta W(\phi) = - \int_{\R^n} \langle H_W,\phi \rangle\,d\mu_W,
$$
\mbox{ where }
$$
H_W(x) = \begin{cases}
H_V(x) & \mbox{ for }x \in \Omega,\\
D\sigma(\sigma(x))H_g(\sigma(x)) &
\mbox{ for }x \in \R^n \backslash \overline{\Omega}\\
\frac{1}{2}\big( H_V(x)  + P H_V(x) + D\sigma(\sigma(x))H_g(\sigma(x))\big)
& \mbox{ for } x \in S.
\end{cases}
$$
This proves that $W$ has weak mean curvature in $L^2(\mu_W)$. Now
$P \nu^S = 0$ on $S$, thus Lemma \ref{lemmaboundarymeancurvature} yields
$$
H_g(x) = H_V(x) - 2 \Delta_P d_S(x) \nu^S(x) \quad \mbox{ $\mu_V$-a.e. on }S.
$$
But $H_W(x) \perp T_x \mu_W$ for $\mu_W$-a.e. $x \in \mathbb{R}^n$ by Brakke's theorem, therefore on $S$ we obtain
$$
0 = P_{T_x \mu_W} H_W(x) = \frac{3}{2} P  H_V(x) \qquad \textrm{for $\mu_W$-a.e. $x \in S$}.
$$
This shows that $H_V(x)$ is perpendicular to $T_x \mu_V$
$\mu_V$-a.e. also on $S$, and claim (b) is proved. Finally claim (c)
follows from Lemma \ref{lemmalocalenergy}. \kasten

%To claim that $B(x,P)(P^\perp v,\,\cdot\,)$ vanishes for all $v \in \R^n$.
%For this we show that $B(x,P)(Pv,\,\cdot\,)$ satisfies the definition of
%weak curvature, then the claim follows by uniqueness, see Proposition
%5.2.2 in \cite{Hut86}. We compute
%\begin{eqnarray*}
%\partial_{P^l_k}\phi^j(x,P) B(x,P)(Pe_j,e_k)^l & = &
%\partial_{P^l_k} P^i_j \phi^j(x,P) B_{ik}^l(x,P)\\
%& = & \partial_{P^l_k} \big(P^i_j \phi^j(x,P)\big) B_{ik}^l
%- \delta_{il} \delta_{jk} \phi^j(x,P)\big) B_{ik}^l\\
%& = & \partial_{P^l_k} (P\phi(x,P))^i B_{ik}^l
%- \phi^k(x,P) \underbrace{B_{lk}^l}_{=0},
%\end{eqnarray*}
%where we used Proposition 5.4(ii) in \cite{Hut86}. Next
%$$
%\langle D_x P \phi, P \rangle
%$$
%$$
%\langle D_x[P^\perp \phi(x)], P \rangle = \langle P^\perp D_x\phi(x),P \rangle = 0.
%$$
%Thus from the definition of varifold curvature, we seem to obtain
%$$
%0 = \int \langle P^\perp {\rm tr\,}B(x,P), \phi(x) \rangle\,dV(x,P).
%$$
%But this is absurd.
%Question Marius: Is it really absurd if eg. our varifold is a minimal surface?
%then $tr B =0$ everywhere, ain't it?

{\em Proof of Lemma \ref{lemmalowerboundary}. } 
We apply the reflection method to $V = V_k$. For $k$ large, 
the assumptions of statement (c) in Proposition \ref{propreflectedvarifold}
are satisfied. Using Willmore's inequality \cite{Wil65}, see
also \cite{KS04}, we get 
\begin{eqnarray*}
16\pi & \leq & \int_{\R^n} |H_{W_k}|^2\,d\mu_{W_k}\\
& \leq & 2 \int_{\overline{\Omega}} |H_{V_k}|^2\,d\mu_{V_k}
+ C \ve \int_{\overline{\Omega}} |B_{V_k}|^2\,d\mu_{V_k}
+ C(\Lambda,\ve) {\bf M}(V_k)\\
& \leq & (2 + C\ve)
\int_{\overline{\Omega}} |B_{V_k}|^2\,d\mu_{V_k}
+ C(\Lambda,\ve) {\bf M}(V_k),
\end{eqnarray*}
 where the last estimate follows from the inequality $\int_{\overline{\Omega}} |H_{V_k}|^2 \, d\mu_{V_k} \leq \int_{\overline{\Omega}} |B_{V_k}|^2 \, d\mu_{V_k}$, which is derived exactly as in (\ref{eqneededbymarius}).
The lemma follows by letting $k \to \infty$ and then
$\ve \searrow 0$.
\kasten

\section{Generic Domains}

In this section we consider orthogonal $(m-1)$-slices $\Gamma$ 
of $S$. By this we mean that $\Gamma$ is a compact, connected, 
$(m-1)$-dimensional submanifold of $S$ along which an affine
$m$-plane intersects $S$ orthogonally. If $\Delta$ is an 
orthogonal $m$-slice of $\Omega$ then each component of 
$\partial \Delta$ has this property, but we do not require
that $\Gamma$ is given like this.

\begin{theorem}[Generic Domains] \label{thmgenericdomains}
Let $\Omega \subset \R^n$ be a bounded $C^2$ domain with
boundary $S = \partial \Omega$. For $m \geq 2$ the 
following holds:
\begin{itemize}
\item[{\rm (1)}] There exist smooth surfaces $S_i$ 
converging to $S$ in $C^2$, which have no orthogonal 
$(m-1)$-slices, in particular the domains $\Omega_i$
have no orthogonal $m$-slices.
\item[{\rm (2)}] The set of domains $\Omega$ without 
orthogonal $m$-slices is open in the $C^2$ topology. 
\end{itemize}
\end{theorem}

Claim (1) does not extend to $m = 1$, in fact for any 
$x_0 \in S$ the line $x_0 + \R \nu^S(x_0)$ intersects
$S$ orthogonally. More interestingly, the segment between
points $x_0, x_1 \in S$ with $|x_0-x_1| = {\rm diam}(S)$
is orthogonal to $S$.

The proof starts with a local analysis near a given 
orthogonal $(m-1)$-slice $\Gamma_0 \subset S$ in the plane 
$z_0+P_0$, where $P_0 \in G(m,n)$, $z_0 \in \R^n$ and $P_0z_0 = 0$.
For any $x \in \Gamma_0$ we have $P_0 = T_x \Gamma_0 \oplus \R\nu^S(x)$
and thus $T_x S = T_x\Gamma_0 \oplus P_0^\perp$. 
Our first goal is to construct a parametrization of all 
nearby $(m-1)$-slices of $S$, whether orthogonal or not.
For suitable $\ve_0 > 0$ we chose Fermi type coordinates
$$
E:\{(x,v) \in \Gamma_0 \times P_0^\perp: |v| < \ve_0\} \to 
U_{\ve_0}(\Gamma_0) \subset S,\,E(x,v) = \pi_S(x+v).
$$
$E$ is a $C^1$ diffeomorphism with $DE(x,0) = {\rm Id}_{T_x S}$.
Describing perturbations of $S$ in direction of the normal $\nu^S$
would lead to a loss of differentiability. Instead we choose
a unit vector field $\xi \in C^2(S,\R^n)$ with 
$\langle \xi,\nu^S \rangle \geq \frac{1}{2}$. For 
$u \in C^2(S)$ we consider the map
$$
g_u \in C^2(S,\R^n),\,g_u(x) = x + u(x)\, \xi(x).
$$
Writing $\nabla$ for the derivative with respect to $x \in S$, we get 
$$
\nabla g_u(x) v  =  v + (du(x)v)  \,  \xi(x) + u(x)  \nabla \xi(x) v.
$$
The set of $u \in C^2(S)$ with $\|u\|_{C^2(S)} < \delta_0$ corresponds 
to a neighborhood of $S$ in the space of surfaces with $C^2$ topology,
sending $u$ to $S_u = g_u(S)$. The unit normal $\nu_u$ along $g_u$
is defined by the equation $N(u,\nu_u) = 0$ where 
$$
N:C^2(S) \times C^1(S,\R^n) \to C^1(TS) \oplus C^1(S),\,
N(u,\nu) = (\nabla g_u)^\ast \nu \oplus \frac{1}{2}(|\nu|^2-1).
$$
Now $\nu_u = \nu^S$ for $u=0$, and we compute 
\begin{eqnarray*}
D_u[\nabla g_u](0)\, \varphi & = & \xi \otimes d\varphi + \varphi\,\nabla \xi,\\
D_u N(0,\nu^S)\, \varphi & = & 
\big(\langle \nu^S,\xi \rangle \,{\rm grad\,}\varphi
+ \varphi (\nabla \xi)^\ast \nu^S\big)\, \oplus \{0\},\\ 
D_\nu N(0,\nu^S)\,\psi & = & P_{TS}\psi \oplus \langle \nu^S, \psi \rangle.
\end{eqnarray*}
Evidently, $D_\nu N(0,\nu^S)$ is an isomorphism, hence 
$\nu:C^2(S) \to C^1(S,\R^n)$, $\nu(u) = \nu_u$, is of class $C^1$. 
Moreover at $u = 0$, its Fr\'{e}chet derivative is 
\begin{equation}
\label{eqvariationnormal}
D\nu(0) \varphi = - \big(\langle \nu^S,\xi \rangle \,{\rm grad\,}\varphi 
+ \varphi (\nabla \xi)^\ast \nu^S\big).
\end{equation}
To describe perturbations of the affine plane, we consider 
the manifold 
$$
M = \{(z,P) \in \R^n \times \R^{n \times n}: P^{{\rm T}} = P = P^2,\,
{\rm tr\,}P = m,\,Pz = 0\}.
$$
%Let $(z(t),P(t)) \in M$ be a curve with $(z(0),P(0)) = (z_0,P_0)$.
%Put $(z'(0),P'(0)) = (\zeta, \Pi)$. By definition of $M$ we then have
%$$
%\Pi^{{\rm T}} = \Pi,\,\,\Pi P_0 + P_0 \Pi = \Pi ,\, \,{\rm tr\,}\Pi = 0,\,\,
%P_0\zeta + \Pi z_0 = 0.
%$$
%The second equation implies $\Pi(P_0) \subset P_0^\perp$ and
%$\Pi(P_0^\perp) \subset P_0$.
%\begin{eqnarray*}
%v \in P_0 & \Rightarrow & P_0 \Pi v = 0 \quad \,\,\Rightarrow \quad
%\Pi v \in P_0^\perp,\\
%v \in P_0^\perp  & \Rightarrow & P_0 \Pi v = \Pi v \quad \Rightarrow \quad
%\Pi v \in P_0.
%\end{eqnarray*}
%Actually, this is an equivalence. 
%Clearly, $\Pi z_0 = - P_0 \zeta$ implies $\zeta + \Pi z_0 = P_0^\perp \zeta$.\\
%\\
%{\bf Step 1: }{\em Local $(m-1)$-slices}\\
For $(z,P,u,v) \in M \times C^2(S) \times C^1(\Gamma_0,P_0^\perp)$ 
with $\|u\|_{C^2(S)} < \delta_0$ and $\|v\|_{C^0(\Gamma_0)} < \ve_0$
we define $F(z,P,u,v) \in C^1(\Gamma_0,P_0^\perp)$ by  
\begin{equation}
\label{eqstep1}
F(z,P,u,v)(x) =
P_0^\perp P^\perp \big(g_u(E(x,v))-z\big)|_{v = v(x)}.
\end{equation}
For $|P-P_0|^2 \leq \frac{1}{2}$ we have 
$|P_0^\perp P^\perp w|^2 \geq \frac{1}{2}\,|P^\perp w|^2$, 
%using $P_0 P^\perp = (P_0-P)P^\perp$, 
which then implies 
$$
F(z,P,u,v) = 0 \quad \Leftrightarrow \quad
g_u\big(E(x,v(x))\big) \in z+P \mbox{ for all }x \in \Gamma_0.
$$
%For the derivative of $F$ with respect to $(z,P) \in M$ at 
%$(z_0,P_0)$ we get 
%\begin{eqnarray*}
%D_{z,P}F(z_0,P_0,0,0)(\zeta,\Pi)(x) & = &
%\frac{d}{dt} F(z(t),P(t),0,0)(x)|_{t=0}\\
%& = & \frac{d}{dt} P_0^\perp P(t)^\perp (x-z(t))|_{t=0}\\
%& = & - P_0^\perp \Pi (x-z_0) - P_0^\perp \zeta.
%%& = & - (\Pi x + \zeta).
%\end{eqnarray*}
%Using $P_0^\perp \Pi (x-z_0) = \Pi(x-z_0)$, since
%$x-z_0 \in P_0$, and $\Pi z_0 + P_0 \zeta = 0$, we can write
%\begin{equation}
%\label{eqFzPderivative}
%D_{z,P}F(z_0,P_0,0,0)(\zeta,\Pi)(x) = - (\Pi x + \zeta) \in P_0^\perp.
%\end{equation}
The function $F$ has derivatives 
\begin{equation}
\label{eqFderivative}
D_u F(z_0,P_0,0,0)\, \varphi = \varphi P_0^\perp \xi 
\quad \mbox{ and }\quad 
D_v F(z_0,P_0,0,0) \eta = \eta.
\end{equation}
%\frac{d}{dt} F(z_0,P_0,t\varphi,0)(x)|_{t=0}
%= \frac{d}{dt} P_0^\perp\big(x+t\varphi(x)\nu^S(x) - z_0\big)|_{t=0}
%Using $\nu^S(x) \in P_0$ we get
%\begin{equation}
%\label{equderivative}
%D_u F(z_0,P_0,0,x,0) \varphi = 0.
%\end{equation}
%Finally for the derivative with respect to $v$ we find, 
%recalling $\eta \in P_0^\perp$,
%$$
%D_v F(z_0,P_0,0,0) \eta(x) =
%\frac{d}{dt} F(z_0,P_0,0,t\eta)|_{t=0} =
%\frac{d}{dt} P_0^\perp \big(E(x,t\eta(x))-z_0\big)|_{t=0} =
%P_0^\perp \eta(x).
%$$
%But $\eta(x) \in P_0^\perp$ by definition, therefore
%\begin{equation}
%\label{eqetaderivative}
%D_v F(z_0,P_0,0,0) \eta = \eta.
%\end{equation}

\begin{lemma} \label{lemmanearbyboundaryslices} There exist 
$\delta_1 \in (0,\delta_0]$, $\ve_1 \in (0,\ve_0]$ and a 
$C^1$ map $\eta:U \to V$ where 
\begin{eqnarray*}
U & = & \{(z,P,u) \in M \times C^2(S): |z-z_0|,\,|P-P_0|,\,
\|u\|_{C^2(S)} < \delta_1\},\\
V & = & \{v \in C^1(\Gamma_0,P_0^\perp): \|v\|_{C^1(\Gamma_0)} < \ve_1\},
\end{eqnarray*}
such that the following statements hold: 
\begin{itemize}
\item[{\rm (1)}] for all $(z,P,u) \in U$, $x \in \Gamma_0$
and $v \in P_0^\perp$ with $|v| < \ve_1$ we have
\begin{equation}
\label{eqlocalslices}
g_u(E(x,v)) \in z+P \quad \Leftrightarrow \quad v = \eta(z,P,u)(x).
\end{equation}
\item[{\rm (2)}] The map $\gamma(z,P,u):\Gamma_0 \to \R^n$, 
$\gamma(z,P,u)(x) = g_u(E(x,v(x)))$ for $v = \eta(z,P,u)$ 
is a $C^1$ parametrization of an $(m-1)$-slice $\Gamma(z,P,u)$
of $S_u$ in the plane $z+P$.
\item[{\rm (3)}] $\eta(z_0,P_0,0) = 0$ and 
$D_u \eta(z_0,P_0,0) \varphi = -\varphi P_0^\perp \xi$.
\end{itemize}
\end{lemma}

{\em Proof. }Applying the implicit function theorem, we obtain 
neighborhoods $U,V$ as claimed such that the equation 
$F(z,P,u,v) = 0$ is solved by a function $v = \eta(z,P,u)$,
and hence $g_u(E(x,v(x))) \in z+P$ for all $x \in \Gamma_0$. Now 
consider
$$
\phi:U \times \Gamma_0 \times P_0^\perp \rightarrow {\mathbb{R}}^n,\,
\phi(z,P,u,x,v) =  P_0^\perp P^\perp\big(g_u(E(x,v))-z\big).
$$
Assume there are sequences $(z_k,P_k,u_k) \to (z_0,P_0,0)$
and $x_k \in \Gamma_0$, such that
$$
\phi(z_k,P_k,u_k,x_k,v_k) = \phi(z_k,P_k,u_k,x_k,w_k) = 0,
$$
where $v_k,w_k \in P_0^\perp$ with $v_k,w_k \to 0$ and $v_k \neq w_k$.
By compactness, we may assume that $x_k \to x_0 \in \Gamma_0$. Then
we obtain by subtracting
\begin{eqnarray*}
0 & = & \phi(z_k,P_k,u_k,x_k,w_k) - \phi(z_k,P_k,u_k,x_k,v_k)\\
& = & D_v \phi(z_0,P_0,0,x_0,0)(w_k-v_k)\\
&& + \int_0^1  \big(D_v \phi(z_k,P_k,u_k,x_k,tw_k +(1-t)v_k)
- D_v \phi(z_0,P_0,0,x_0,0)\big) \cdot (w_k-v_k)\,dt\\
& = & w_k-v_k + o(|w_k-v_k|).
\end{eqnarray*}
For $k \to \infty$ we get a contradiction, and the
lemma follows. \kasten

%Differentiating the equation $F(z,P,u,\eta(z,P,u)) = 0$
%we obtain
%\begin{eqnarray}
%\label{eqetazPderivative}
%D_{z,P} \eta(z_0,P_0,0) (\zeta,\Pi)(x) & = & \Pi x+\zeta,\\
%\begin{equation}
%\label{eqetauderivative}
%D_u \eta(z_0,P_0,0)(x) \varphi = -\varphi P_0^\perp \xi.
%\end{equation}

%D_x \eta(z_0,P_0,0)(x) & = & 0.
%\end{eqnarray}
%Now let $\gamma(t)$ be a curve in $\Gamma_0$ such that $\gamma(0) = x$
%and $\gamma'(0) = \xi$. Then we have
%$$
%D_x F(z_0,P_0,0,x,0) \xi = \frac{d}{dt}P_0^\perp(\gamma(t)-z_0)_{t=0}
%= P_0^\perp \xi.
%$$
%But since $\gamma(t) \in z_0+ P_0$, we have $\xi \in P_0$ and hence
%\begin{equation}
%\label{eqxderivative}
%D_x F(z_0,P_0,0,x,0) \xi = 0.
%\end{equation}
%In fact differentiating the equation $F(z,P,u,\eta(z,P,u)) = 0$ we have
%\begin{eqnarray*}
%0 & = & D_{z,P} F(z_0,P_0,0,0) (\eta,\Pi)(x) +
%D_v F(z_0,P_0,0,0) D_{z,P} \eta(z_0,P_0,0)(x)\\
%& = & - (\Pi x + \zeta) + D_{z,P} \eta(z_0,P_0,0) (\zeta,\pi)(x).
%\end{eqnarray*}
%\begin{eqnarray*}
%0 & = & D_u  F(z_0,P_0,0,x,0) \varphi +
%D_v F(z_0,P_0,0,0) D_u \eta(z_0,P_0,0) \varphi\\
%& = & D_u \eta(z_0,P_0,0,x) \varphi.
%\end{eqnarray*}
%\begin{eqnarray*}
%0 & = & D_x F(z_0,P_0,0,x,0) \xi +
%D_v F(z_0,P_0,0,x,0) D_x D_x \eta(z_0,P_0,0,x) \xi\\
%& = & D_x \eta(z_0,P_0,0,x) \xi.\\
%\end{eqnarray*}
We now turn to orthogonal slices near $\Gamma_0$. For $U$ 
as in Lemma \ref{lemmanearbyboundaryslices} we define 
the map $G:U \to C^1(\Gamma_0,P_0^\perp)$ by 
\begin{equation}
G(z,P,u)(x) = (P_0^\perp P^\perp \nu_u)(E(x,v))|_{v = \eta(z,P,u)(x)}.
\end{equation}
We have $G(z,P,u) = 0$ if and only if $\gamma(z,P,u)$ 
parametrizes an orthogonal $(m-1)$-slice of $S_u$. 
%in particular $G(z_0,P_0,0)(x) = P_0^\perp \nu^S(x) = 0$ by assumption.
%Let us compute the derivatives of $G$. For $(z(t),P(t)) \in M$
%with $(z(0),P(0)) = (z_0,P_0)$, $(z'(0),P'(0)) = (\zeta,\Pi)$
%we have
%\begin{eqnarray*}
%D_{z,P}G(z_0,P_0,0)(\zeta,\Pi)(x) & = &
%\frac{d}{dt} G(z(t),P(t),0)(x)|_{t=0}\\
%& = & \frac{d}{dt}P_0^\perp P(t)^\perp \nu^S(E(x,v))|_{v = \eta(z(t),P(t),0)(x)}\\
%& = & - P_0^\perp \Pi\, \nu^S(x) + P_0^\perp D\nu^S(x) D\eta(z_0,P_0,0)(\zeta,\Pi)(x)\\
%& = & - \Pi\, \nu^S(x) + P_0^\perp D\nu^S(x) (\Pi x + \zeta).
%\end{eqnarray*}
%For a curve $\gamma:(-\ve,\ve) \to \Gamma_0$ with $\gamma(0) = x$ we have
%$$
%0 = P_0^\perp \nu^S(\gamma(t))|_{t=0} = P_0^\perp D\nu^S(x) \gamma'(0).
%$$
%As $D\nu(x)\gamma'(0) \perp \nu^S(x)$ we see
%$D\nu^S(x)T_x \Gamma_0 \subset T_x \Gamma_0$. Now $D\nu^S(x)$ is
%symmetric and thus $D\nu^S(x)P_0^\perp \subset P_0^\perp$. We
%conclude for $x \in \Gamma_0$
%\begin{equation}
%\label{eqGzPderivative}
%D_{z,P}G(z_0,P_0,0)(\zeta,\Pi) =
%-\Pi\,\nu^S(x) + D\nu^S(x) (\Pi x + \zeta) \in P_0^\perp.
%\quad \mbox{ for }x \in \Gamma_0.
%\end{equation}
The derivative $D_u G(z_0,P_0,0)$ is the operator 
$L_0:C^2(S) \to C^1(\Gamma_0,P_0^\perp)$ given by 
\begin{equation}
\label{eqGuderivative}
L_0\varphi = 
- P_0^\perp \Big(\langle \nu^S,\xi \rangle\,{\rm grad\,}\varphi 
+ \varphi \,(\nabla \xi)^\ast \nu^S
+ \varphi\, \xi\Big)\big\vert_{\Gamma_0}. 
\end{equation}
%& = & \frac{d}{dt}\, P_0^\perp
%\nu(t \varphi)(E(x,v)|_{v = \eta(z_0,P_0,t\varphi)(x)}\\
%& = & - P_0^\perp \big(\langle \nu^S,\xi \rangle\, {\rm grad\,}\varphi 
%+ \varphi\, (\nabla \xi)^\ast \nu^S\big)\\ 
%&& + P_0^\perp D\nu^S(x) D_u \eta(z_0,P_0,0) \varphi\\
%As $D_u \eta(z_0,P_0,0) = - P_0^\perp \xi$ we arrive at
%\begin{equation}
%\label{eqGuderivative}
%D_u G(z_0,P_0,0) \varphi = P_0^\perp {\rm grad\,} \varphi|_{\Gamma_0}.
%\end{equation}
In Appendix B we construct a right inverse 
$R_0:C^1(\Gamma_0,P_0^\perp) \to C^2(S)$ to $L_0$. 
The space $W_0 = {\rm im\,}R_0$ is closed in $C^2(S)$, 
and we have the direct sum decomposition
\begin{equation}\label{expliziesKOmplementvonKerL0}
C^2(S) = {\rm ker\,}L_0 \oplus W_0,\,
\varphi = \varphi - R_0(L_0\varphi) \oplus R_0(L_0\varphi).
\end{equation}
Thus $L_0|_{W_0}:W_0 \to C^1(\Gamma_0,P_0^\perp)$
is an isomorphism. 

\begin{lemma} \label{lemmalocalexcept}
There exist $\delta_2 \in (0,\delta_1]$, $\ve_2 \in (0,\ve_1]$ and a 
map $g_0 \in C^1(A_0,W_0)$ where 
$$
A_0 = \{(z,P,\chi) \in M \times {\rm ker\,}L_0:
|z-z_0|,\,|P-P_0|,\,\|\chi\|_{C^2(S)} < \delta_2\},
%B_0 & = & \{w \in W_0: \|w\|_{C^2(S)} < \ve_2\},
$$
such that for any $u = \chi + w \in {\rm ker\,}L_0 \oplus W_0$ 
with $\|u\|_{C^2(S)} < \ve_2$ we have the equivalence
$$
G(z,P,u) = 0 \quad \Leftrightarrow \quad
\|\chi\|_{C^2(S)} < \delta_2 \mbox{ and }w = g_0(z,P,\chi). 
$$
\end{lemma}

%By assumption $g(z_0,P_0,0) = 0$. For the $(z,P)$ derivative we have 
%\begin{eqnarray*}
%0 & = & \frac{d}{dt} G(z(t),P(t),g(z(t),P(t),0))|_{t=0} \\
%& = & D_{z,P}G(z_0,P_0,0) (\zeta,\Pi)
%+ D_u G(z_0,P_0,0) D_{z,P}g(z_0,P_0,0) (\zeta,\Pi)\\
%& = & -\Pi\,\nu^S + D\nu^S (\Pi\, {\rm id}_S + \zeta)
%+ L_0\,D_{z,P}g(z_0,P_0,0) (\zeta,\Pi).
%\end{eqnarray*}
%Thus $\varphi = D_{z,P}g(z_0,P_0,0) (\zeta,\Pi)$ is the unique
%element of $W_0$ with 
%\begin{equation}
%\label{eqgzPderivative}
%L_0 \varphi = \big(\Pi\,\nu^S - D\nu^S (\Pi\, {\rm id}_S + \zeta\big)|_{\Gamma_0}.
%\end{equation}
%The derivative in the direction of $\chi \in {\rm ker\,}L_0$ is 
%$D_\chi g(z_0,P_0,0) = 0$, in particular we can assume by 
%taking $\delta_0 > 0$ slightly smaller that
%$$
%0 = \frac{d}{dt} G\big(z_0,P_0,t\chi + g(z_0,P_0,t\chi)\big)|_{t=0}
%= L(\chi + D_v g(z_0,P_0)\chi).
%$$
%But $L \chi = 0$ and $L$ is injective on $W$, therefore 
%\begin{equation}
%\label{eqGxiderivative}
%D_\chi g(z_0,P_0,0) = 0.
%\end{equation}
We define the family of functions 
\begin{equation}
\label{eqfamily}
C(\Gamma_0) = \{u = \chi + g_0(z,P,\chi): 
|z-z_0|,|P-P_0| \leq  \frac{\delta_2}{2},\,\|\chi\|_{C^2(S)} < \delta_2\}.
\end{equation}
For $\|u\| < \ve_2$ we have $u \in C(\Gamma_0)$ if and only if the
slice $\Gamma(z,P,u)$ is orthogonal to $S_u$, for some 
$(z,P) \in M$ with $|z-z_0|,|P-P_0| \leq \frac{\delta_2}{2}$. 
%For $\phi(z,P,\chi) = \chi + g(z,P,\chi)$ we have $\phi(z_0,P_0,0) = 0$
%and $\phi$ has derivatives
%\begin{eqnarray}
%\label{eqphizPderivative}
%D_{z,P}\phi(z_0,P_0,0) & = & D_{z,P}g(z_0,P_0,0),\\
%\label{eqphivderivative}
%D_\chi \phi(z_0,P_0,0) & = & {\rm Id}_{{\rm ker\,}L_0}. 
%\end{eqnarray}

In the following we employ the Hausdorff distance 
$d_H(\Gamma,\Gamma')$ for slices $\Gamma,\Gamma'$.

\begin{lemma} \label{lemmahausdorff}
There is a constant $\delta = \delta(\Gamma_0) > 0$ such that 
whenever $\Gamma$ is an orthogonal $(m-1)$-slice of $S_u$ in 
$z+P$ with $d_H(\Gamma,\Gamma_0),\,\|u\|_{C^2(S)} < \delta$ 
then $u \in C(\Gamma_0)$.
\end{lemma}

{\em Proof. }Let $\Gamma_k \subset z_k + P_k$ be orthogonal  $(m-1)$-slices
of $S_{u_k}$, where $\mathrm{dist}_H(\Gamma_0,\Gamma_k) \to 0$
and $\|u_k\|_{C^2(S)} \to 0$. Then $(z_k,P_k)$ (sub-)converges
to some $(z,P) \in M$, where $\Gamma_0 \subset z+P$ by
Hausdorff convergence. We must have $(z,P) = (z_0,P_0)$,
otherwise $\Gamma_0$ would be contained in an affine subspace 
of dimension at most $(m-1)$. For $k$ large we have 
$\Gamma_k \subset U_{\ve_1}(\Gamma_0)$. Lemma  
\ref{lemmanearbyboundaryslices}(1) now implies that 
$\Gamma_k = \Gamma(z_k,P_k,u_k)$. We have 
$\eta(z_k,P_k,u_k) \to 0$ in $C^1(\Gamma_0,P_0^\perp)$.
Now decompose $u_k = \chi_k + w_k \in {\rm ker\,}L_0 \oplus W_0$.
Then $|z_k-z_0|,|P_k-P_0| \leq \frac{\delta_2}{2}$, 
$\|u_k\|_{C^2(S)} < \ve_2$ and $\|\chi_k\|_{C^2(S)} < \delta_2$ 
for $k$ large. As $\Gamma_k$ is an orthogonal slice, we have
$$
G(z_k,P_k,u_k)(x) = 
(P_0^\perp P^\perp \nu_{u_k})(E(x,v))|_{v = \eta(z_k,P_k,u_k)(x)}
= 0 \quad \mbox{ for all }x \in \Gamma_0.
$$
Therefore Lemma \ref{lemmalocalexcept} and (\ref{eqfamily}) 
yield $u_k \in C(\Gamma_0)$. \kasten

This finishes the local analysis near $\Gamma_0$. 
The following lemma implies that the set ${\cal O}(S)$
of orthogonal $(m-1)$-slices of $S$, equipped with 
Hausdorff distance, becomes a compact metric space.

\begin{lemma} \label{lemmahausdorffcompact} Let $\Gamma_k$ be orthogonal 
$(m-1)$-slices for $S_{u_k}$, where $\|u_k\|_{C^2(S)} \to 0$. 
After passing to a subsequence, the $\Gamma_k$ converge to an
orthogonal $(m-1)$-slice $\Gamma_0 \neq \emptyset$ of $S$
with respect to Hausdorff distance.
\end{lemma}

{\em Proof. }We assume that $\Gamma_k$ is contained in $z_k+P_k$
where $P_kz_k = 0$. Arguing as in the proof of Theorem
\ref{thmboundsnearboundary}, we get a point $x_k \in z_k+P_k$
in $\Omega_{u_k}$ with $$\mathrm{dist}(x_k,S_{u_k}) \geq \varrho_{S_{u_k}} \xrightarrow{k \rightarrow \infty} \rho_S >0.$$ We may
assume $P_k \to P_0$ and $x_k \to a \in \Omega$, which implies
$z_k = P_k^\perp x_k \to P_0^\perp a =:z_0$. Now consider a point
$y_0 = \lim_{k \to \infty} y_k \in S$ where $y_k \in \Gamma_k$.
Putting $S_k = S_{u_k}$ we have 
$\nu^S(y_0) = \lim_{k \to \infty} \nu^{S_k}(y_k) \in z_0+P_0$.
Let $A^{\Gamma_k}$ be the second fundamental form of 
$\Gamma_k$ as submanifold of $\R^n$. As $\Gamma_k$ is totally 
geodesic in $S_k$, we have
$$
A^{\Gamma_k}(y)(v,w) = A^{S_k}(v,w) \quad \mbox{ for any }v,w \in T_y\Gamma_k.
$$
Hence the norms $\|A^{\Gamma_k}\|_{L^\infty}$ are bounded.
Moreover, Lemma \ref{lemmamasses} implies a bound for the 
measures $|\Gamma_k|$. After passing to a subsequence, we can assume 
$$
{\cal H}^{m-1} \llcorner \Gamma_k \to \gamma \quad
\mbox{ in }C^0_c(\R^n)'.
$$
We claim that $\Gamma_k$ converges to $\Gamma:={\rm spt\,}\gamma$
in Hausdorff distance. First assume by contradiction that there
exist $x_k \in \Gamma$ with
$$
\limsup_{k \to \infty} \mathrm{dist}(x_k,\Gamma_k) \geq \delta > 0.
$$
After passing to a subsequence, we get $x_k \to x \in \Gamma$
and $\mathrm{dist}(x,\Gamma_k) \geq \delta > 0$, a contradiction.
Secondly assume, again by contradiction, that there are
$x_k \in \Gamma_k$ with
$$
\limsup_{k \to \infty} \mathrm{dist}(x_k,\Gamma) \geq \delta > 0.
$$
There is a $\varrho > 0$ (the Langer radius) such that
for any $k$ and any $x \in \Gamma_k$, the submanifold $\Gamma_k$
contains a graphical piece given by a function
$$
h_k:T_x \Gamma_k \cap B_\varrho(x) \to T_x\Gamma_k^\perp
\quad \mbox{ where }|Dh_k| \leq 1,
$$
see for instance Theorem 2.6 in \cite{Bre15}. In particular
$|\Gamma_k \cap B_{\frac{\delta}{2}}(x_k)| \geq c > 0$
uniformly, again a contradiction. It remains to show that 
$\Gamma \subset S \cap (z_0+P_0)$ is an $(m-1)$-dimensional 
manifold. This follows by representing the slices $\Gamma_k$ 
locally as graphs, using the implicit function theorem as before. 
\kasten

{\em Proof of Theorem \ref{thmgenericdomains}.} By Lemma \ref{lemmahausdorffcompact} the set ${\cal O}(S)$ has a finite 
covering by the Hausdorff distance balls
$$
{\cal O}(S) = \bigcup_{i=1}^\ell B_{\delta_i}(\Gamma_i)
\quad \mbox{ for }\delta_i = \delta(\Gamma_i) > 0 
\mbox{ as in Lemma } \ref{lemmahausdorff}.
$$
We claim that if $S_u$ has an orthogonal $(m-1)$-slice 
$\Gamma$ for $\|u\|_{C^2(S)}$ small, then 
$$
u \in \bigcup_{i=1}^\ell C_i \quad \mbox{ where }C_i = C(\Gamma_i).
$$
Namely, let $\Gamma'_k$ be orthogonal $(m-1)$-slices 
of $S_{u_k}$ where $u_k \to 0$ in $C^2(S)$. By 
Lemma \ref{lemmahausdorffcompact} we can assume that 
the $\Gamma'_k$ converge to an orthogonal $(m-1)$-slice
$\Gamma$ of $S$. Now $\Gamma \in B_{\delta_i}(\Gamma_i)$
for some $i \in \{1,\ldots,\ell\}$. It follows that
$d_H(\Gamma'_k,\Gamma_i) < \delta_i$ for $k$ large.
But then $u_k \in C_i$ by Lemma 
\ref{lemmahausdorff}.\\
\\
%Now for $i = 1,\ldots,\ell$ let $\delta_i,\ve_i > 0$ and 
%$g_i \in C^1(A_i,B_i)$ as in Lemma \ref{lemmalocalexcept}. 
%Put
%$$
%{\cal F}_i = {\cal F}(\Gamma_i) = 
%\{u = \chi + g_i(z,P,\chi): (z,P,\chi) \in A_i\}.
%$$
It is easy to see that $C_i \cap B_\varrho(0)$ is closed in 
$B_\varrho(0)$ for $\varrho > 0$ small. 
%For this assume that $u_k \in C_i$ 
%converges to $u \in B_\varrho(0)$ where $u_k = \chi_k + w_k$ 
%with $w_k = g_i(z_k,P_k,\chi_k)$. Then $\chi_k \to \chi$ 
%and $w_k \to w$. 
%$(z_k,P_k) \to (z,P) \in M_i$. But then
%$$
%w = \lim_{k \to \infty} g_i(z_k,P_k,\chi_k) = g(z,P,\chi).
%$$
For given $u \in B_\varrho(0) \cap C_i$ we now construct 
a perturbation $u' \in U_i = B_\varrho(0) \backslash C_i$. 
Decompose $u = \chi + w$ and choose a subspace $X \subset W_i$ of 
dimension $\dim M + 1$ with $w \in X$. Here $W_i$ is the space introduced in \eqref{expliziesKOmplementvonKerL0} for $\Gamma_i$. By Hahn-Banach
there is a continuous projection $Q :W_i \to X$. Consider
the $C^1$ map
$$
\phi: M' = \big\{(z,P) \in M: |z-z_0|,|P-P_0| \leq \frac{\delta_2}{2}\big\} \to X,\,
\phi(z,P) = Q(g_i(z,P,\chi)). 
$$
The image of $\phi$ is a Lebesgue null set in $X$. Thus 
there exists $w' \in X \cap B_\ve(w)$ with 
$$
w' \,\notin \, Q(g_i(M',\chi)), 
\mbox{ in particular } 
w' \notin g_i(M',\chi) \cap X.
$$
It follows that $u' = \chi + w' \in U_i$. Moreover, for 
those $j \in \{1,\ldots,\ell\}$ for which $u$ is already 
in $U_j$, we can arrange that also $u' \in U_j$ by choosing 
the perturbation sufficiently small. Repeating this step at
most $\ell$ times yields a function $u_\varrho \in \bigcap_{i=1}^\ell U_i$, 
which means that $S_{u_\varrho}$ has no orthogonal slices. 
By approximation we even obtain a smooth surface with that 
property, applying Lemma \ref{lemmahausdorffcompact}. 
This finishes the proof of statement (1) in Theorem 
\ref{thmgenericdomains}.\\
\\
For statement (2) consider a sequence $\Omega_{u_k}$ with 
$\|u_k\|_{C^2(S)} \to 0$, where $\Omega_{u_k}$ has an 
orthogonal $m$-slice $\Delta_k$ in $z_k+P_k$. By Lemma 
\ref{lemmahausdorffcompact} we may assume 
$(z_k,P_k) \to (z,P)$, and $\Gamma_k  = \partial \Delta_k$
converges in Hausdorff distance to an orthogonal $(m-1)$-slice 
$\Gamma$ (possibly with several components) of $S$. Now put
$$
\Delta = \{x \in (z+P) \backslash \Gamma: x = \lim_{k \to \infty} x_k
\mbox{ for some }x_k \in \Delta_k\}. 
$$
By the argument of Theorem \ref{thmboundsnearboundary} which was
also used in Lemma \ref{lemmahausdorffcompact}, there are points $a_k \in \Delta_k$ 
with $\mathrm{dist}(a_k,S) \geq \varrho_S > 0$, thus $\Delta \neq \emptyset$.
Now if $x \in \Delta$ and $\varrho < \mathrm{dist}(x,\Gamma)$, then  
$(z_k+P_k) \cap B_\varrho(x) \subset \Delta_k$ for $k$ large, 
hence $(z+P) \cap B_\varrho(x) \subset \Delta$. This shows that
$\Delta$ is open and relatively closed in $(z+P) \backslash \Gamma$, 
and hence $\Delta$ is an orthogonal $m$-slice of $\Omega$. The
proof of Theorem \ref{thmgenericdomains} is now complete.

\section{Appendix}

{\bf Appendix A. }The following disintegration theorem is used 
in the paper. We refer to \cite{Sim84}, Lemma 38.4, and \cite{Amb05}, 
Theorem 5.3.1, for further information. 

\begin{theorem}[disintegration] \label{lemmadisintegration} 
Let $\gamma$ be a Radon measure on $X \times Y$ with compact
support, where $X,Y$ are metric spaces, and denote by 
$\pi:X \times Y \to Y$ the projection map. Then
$\beta = \pi_\ast \gamma$ is a Radon measure, and 
there is a $\beta$-a.e. uniquely determined family 
$(\alpha_y)_{y \in Y}$ of Radon probability measures 
on $X$, such that for any Borel function
$\phi:X \times Y \to [0,\infty]$  one has 
\begin{equation}
\label{eqdisintegration}
\int_{X \times Y} \phi\,d\gamma = 
\int_Y \int_X \phi(x,y)\,d\alpha_y(x)\,d\beta(y).
\end{equation} 
\end{theorem}

{\bf Appendix B. }In this appendix, we construct the (bounded, linear) right inverse $R_0$ 
as announced in equation \eqref{eqGuderivative}. We first construct $R_0$ on $\R^n = \R^m \times \R^p$. 

\begin{lemma} \label{lemmaextension} There exists a bounded operator 
$R_0:C^0(\R^m,\R^p) \to C^1(\R^n)$ such that 
$$
R_0 \psi(x,0) = 0 \quad \mbox{ and } \quad 
{\rm grad\,}(R_0 \psi)(x,0) = (0,\psi(x)) \quad \mbox{ for all }x \in \R^m. 
$$
Moreover, $\|R_0 \psi\|_{C^{k+1}(\R^n)} \leq C_k\,\|\psi\|_{C^k(\R^m)}$ 
and $(R_0\psi)(x,z) = 0$ for $\mathrm{dist}(x,{\rm spt\,}\psi) \geq |z|$.
\end{lemma}

{\em Proof. } Fix $\eta \in C^1_c(\R^m)$ with support in $B_1(0)$ and 
weight $\int_{\R^m} \eta(x)\,dx = 1$. For $\psi  \in C^0(\R^m)$ we define
$u \in C^0(\R^m \times \R)$ by 
$$
u(x,z) = z\,(\eta_z \ast \psi)(x) \quad \mbox{ where }
\eta_z \ast \psi(x) = \int_{\R^m} \eta(y) \psi(x-zy)\,dy.
$$
Clearly $u(x,0) = 0$. For $z \neq 0$  we have, putting 
$\sigma_{m,z} = ({\rm sign\,} z)^m$, 
%\begin{equation}
%\label{eqapp1}
%|u(x,z)| \leq C |z|\, \|\psi\|_{C^0(B_{|z|}(x))} \leq C|z|\, \|\psi\|_{C^0(\R^m)}. 
%\end{equation}
$$
(\eta_z \ast \psi)(x) = \sigma_{m,z} 
\int_{\R^m} z^{-m} \eta\Big(\frac{x-y}{z}\Big)\,\psi(y)\,dy.
$$
Differentiating under the integral and transforming back yields 
\begin{eqnarray*}
\partial_i u(x,z) & = & \int_{\R^m} \partial_i \eta(y)\,\psi(x-z y)\,dy
\quad \mbox{ for }i = 1,\ldots,m,\\
%\sigma_{m,z_i} \int_{\R^m}
%z_i^{-m} D\eta\Big(\frac{x-y}{z_i}\Big)\,\psi_i(y)\,dy
\partial_z u(x,z) & = & 
%\eta_{z_i} \ast \psi_i(x) - \int_{\R^m} z_i^{-m}
%\Big(m \eta\Big(\frac{x-y}{z_i}\Big) 
%+ D\eta\Big(\frac{x-y}{z_i}\Big) \cdot \frac{x-y}{z_i}\Big)
%\psi(y)\,dy\\
\int_{\R^m} \eta(y) \psi(x-z y)\,dy
- \int_{\R^m} \Big(m \eta(y) + \sum_{i=1}^m y^i \partial_i \eta(y)\Big)\,\psi(x-z y)\,dy.
\end{eqnarray*}
Thus $Du \in C^0(\R^{m+1},\R^{m+1})$ and 
$\|u\|_{C^1(\R^{m+1})} \leq C\,\|\psi\|_{C^0(\R^m)}$. Since $\eta$ 
has weight one while $m \eta + y \cdot D\eta = \mathrm{div} \big(\eta(y)y\big)$
has weight zero, we conclude $\partial_z u(x,0) = \psi(x)$. Now 
assume $\psi \in C^k(\R^m)$ for some $k \geq 1$. We have 
for any $\alpha \in \N_0^m$, $s \in \N_0$ with $|\alpha| + s = k$
\begin{eqnarray*}
\partial^\alpha (\partial_z)^s  (\eta_z \ast \psi)(x,z) & = & 
\int_{\R^m} \sum_{\beta \in \N_0^m: |\beta| = s} 
\eta_\beta(y)\,\partial^{\alpha + \beta} \psi(x-zy)\,dy, 
\quad \mbox{ where }\\
\eta_\beta(y) & = & (-1)^{|\beta|} \frac{s!}{\beta!} \eta(y) y^\beta.
\end{eqnarray*}
For $s = 0$ we have $\partial^\alpha (z\, \eta_z \ast \psi) 
= z\,\partial^\alpha (\eta_z \ast \psi)
= z\,\eta_z \ast \partial^\alpha \psi$,
while for $s \geq 1$ 
\begin{eqnarray*}
\partial^\alpha (\partial_z)^s (z\, \eta_z \ast \psi) & = & 
z\, \partial^\alpha (\partial_z)^s (\eta_z \ast \psi)
+ s \, \partial^\alpha (\partial_z)^{s-1} (\eta_z \ast \psi)\\
& = & \sum_{\beta \in \N_0^m: |\beta| = s} 
z\,(\eta_\beta)_z \ast \partial^{\alpha + \beta}\psi 
+ s  \sum_{\beta \in \N_0^m:|\beta| = s-1} 
(\eta_\beta)_z \ast \partial^{\alpha + \beta}\psi.
\end{eqnarray*}
The second sum involves only derivatives of order at most $k-1$. For 
the first sum we apply the argument above, with $\eta$ replaced
by $\eta_\beta$. This shows 
$\partial^\alpha (\partial_z)^s (z \,\eta_z \ast \psi) \in C^1(\R^{m+1})$
and hence $z \,\eta_z \ast \psi \in C^{k+1}(\R^{m+1})$.
To prove the theorem, we finally put $u(x,z) = \sum_{j=1}^p u_j(x,z)$ where 
$u_j(x,z) = \sum_{j=1}^p z_j\,\eta_{z_j} \ast \psi_j(x)$. 
The claim follows immediately defining $R_0 \psi = u$. \kasten 

Before we proceed to the general case we observe that if ${\rm spt}(\psi_j) \subset W$ for all $j=1,...,p$ then for each $\epsilon > 0$ one may also achieve that ${ \rm spt  }(R_0\psi) \subset B_{2\sqrt{p}\epsilon}(W)\times (-\epsilon,\epsilon)^p$. 
%Indeed, let $\alpha \in C_0^\infty((-\epsilon,\epsilon)^p)$ be a fixed function such that $\alpha \equiv 1$ on $(- \frac{\epsilon}{2},\frac{\epsilon}{2})^p$. 
 Indeed, let $\alpha \in C_0^\infty((-\epsilon,\epsilon)^p)$ be a fixed function such that $\alpha \equiv 1$ on $(- \frac{\epsilon}{2},\frac{\epsilon}{2})^p$ and look at $\tilde{R}_0 : C^0(\R^m, \R^p) \rightarrow C^1(\R^n)$ defined via 
\begin{equation}\label{eq:R0schlange}
(\tilde{R}_0 \psi )(x,z) = \alpha(z) (R_0 \psi)(x,z). 
\end{equation}
Clearly, $\tilde{R}_0$ has all the properties in the statement of the previous lemma. Moreover, $(x,z) \in {\rm spt}( \tilde{R}_0 \psi)$ implies $z \in (-\epsilon,\epsilon)^p$ and $(x,z) \in {\rm spt}(R_0 \psi)$. One infers $|z| \leq \sqrt{p}\,\epsilon$  and by the previous lemma this means ${\rm dist}(x, {\rm spt}(\psi) ) \leq \sqrt{p}\,\epsilon$.

 %Indeed, instead of our choice in the above proof one could repeat the construction of $R_0$ with 
%\begin{equation*}
%u_j(x,z) = \sum_{i=1}^p  \alpha(z) z (\eta_{z} * \psi_i)(x), 
%\end{equation*}
%where $\alpha \in C_c^\infty((-\epsilon,\epsilon))$ is a fixed function such that $\alpha \equiv 1$ on $(- \frac{\epsilon}{2}, \frac{\epsilon}{2})$. 
%Now we turn to the general situation, which finally constructs $R_0$ in the situation of \eqref{eqGuderivative}.

\begin{lemma}\label{lem:7.1}
Let $N$ be an $n$-dimensional $C^l$-manifold and $M$ be a compact $m$-dimensional $C^l$-submanifold without boundary, $l \geq 1$. Then there exists a bounded linear operator $R_0 : C^0(M, (TM)^\perp) \rightarrow C^1(N)$ such that 
\begin{equation*}
R_0 \psi \vert_M = 0 \quad {\rm and} \quad {\rm grad\,} R_0 \psi \vert_M = \psi.
\end{equation*}
Further $\|R_0 \psi\|_{C^{k+1}(N)} \leq C_k(M,N) \|\psi\|_{C^k(M)}$ for all $k = 0,...,l-1$. 
%Further, if $\psi \in C^{k}(M, (TM)^\perp)$ for any $k =1,...,l-1$ it is possible to choose $u(\psi) \in C^{k+1}(N)$. Moreover, one can achieve that the map $T: C^k(M, (TM)^\perp) \rightarrow C^{k+1}(N), \psi \mapsto u(\psi)$ defines a continuous linear operator.
\end{lemma}

{\it Proof of Lemma \ref{lem:7.1}.}
Let $M$ and $N$ be as in the statement and $\psi \in C^{k-1}(M,(TM)^\perp)$ be arbitrary.  
 It is possible to choose a finite collection of open sets $U_1,...,U_r \subset N$ covering $M$ and maps $\Phi_j \in C^l(U_j;\R^n)$ which are diffeomorphisms onto their images and satisfy $\Phi_j(U_j \cap M) = \Phi_j(U_j) \cap \R^m$.  Note in particular that $D\Phi_j(x)(T_xM) = \R^m$. 
% and $D\Phi_j(x)(T_xM^\perp)= \{ 0 \} \times \R^p$ for all $x \in U_j \cap M$.
 Let now $\beta_1,...,\beta_r$ be a partition of unity for $M$ such that $\beta_j \in C^l_c(M \cap U_j)$ for all $j$. By openness of $U_j$ it is possible to choose $\epsilon > 0$ such that $B_{2\sqrt{p}\epsilon}(\Phi_j({ \rm spt } ( \beta_j ) )) \times (-\epsilon,\epsilon)^p \subset \subset \Phi_j(U_j)$ for all $j =1,...,r$. Next define $w_j := \beta_j [(D\Phi_j)^*]^{-1} \psi \circ \Phi_j^{-1}$ and
 $u : N \rightarrow \R$ to be 
\begin{equation*}
u(x)  = \sum_{j = 1}^r \tilde{R}_0 w_j (\Phi_j(x)) ,
\end{equation*}
where $\tilde{R}_0$ is chosen as in \eqref{eq:R0schlange} with $\epsilon> 0$ as above, in particular one has ${ \rm spt }(\tilde{R}_0 w_j)\subset \subset \Phi_j(U_j)$.
We remark that for the well-definedness  and $C^{k+1}$-smoothness of $\tilde{R}_0 w_j$ we need to ensure that  $w_j \in C^k(\R^m,\R^p)$. This is due to the fact that $D \Phi_j \in C^{l-1}$ and maps $TM$ to $\R^m$, which implies that $(D\Phi_j)^*$ maps $\R^p$ to $(TM)^\perp$.
 One readily checks that $u \vert_M = 0$ and with the chain rule one obtains for $x \in M$  
\begin{align*}
{\rm grad }(u)(x)  & = \sum_{j = 1}^r D\Phi_j(x)^* ( {\rm grad} (\tilde{R}_0 w_j)  (\Phi_j(x)) \\  &= \sum_{j = 1}^r D\Phi_j(x)^* (  \beta_j(x)  (D\Phi_j(x)^*)^{-1} \psi(x)  ) 
= \psi(x).  
\end{align*}
The claim follows now defining $R_0 \psi = u$. From the explicit construction one 
easily deduces the desired operator norm estimates for $R_0$. \kasten

%Now notice that for each $x \in M$  $D\Phi_j(x)^*$ maps $\R^m \times \{0 \}$ to $T_xM$, since $D\Phi_j(x)$ maps $(T_xM)^\perp$ to $\{ 0 \} \times \R^p$. Since $R_j(\Phi_j(x)) \in \R^m \times \{ 0 \} $ for all $x \in M$ one infers from the previous equation
%In particular for all $x \in M$ one has
%\begin{equation*}
%\Pi_{(TM)^\perp} {\rm grad }(u)(x) = \Pi_{(TM)^\perp}  \psi(x) = \psi (x).  
%\end{equation*}
%Linearity and conitnuity of the associated operator follows directly from our explicit construction. \kasten

\vspace{1cm} Mathematisches Institut der\\Albert-Ludwigs-Universit\"at Freiburg\\Ernst-Zermelo-Stra\ss e 1,
D-79104 Freiburg,
Germany\\email: ernst.kuwert@math.uni-freiburg.de. \vspace{1cm}\\
Institut für Mathematik der\\ Universität Augsburg \\ Universitätsstra\ss e 14,
D-86159 Augsburg,
Germany\\email: marius1.mueller@uni-a.de. 

\end{document}